\documentclass{article}
\usepackage{graphicx}
\usepackage{amsmath}
\usepackage{amsfonts}
\usepackage{amssymb}

\newtheorem{proposition}{Proposition}  

\newenvironment{proof}[1][Proof]{\textbf{#1.} }{\ \rule{0.5em}{0.5em}}

\begin{document}

\title{Ricci flow on locally homogeneous closed 4-manifolds}
\date{}
\author{
James Isenberg \\
{\small Department of Mathematics and Institute of Theoretical
Science}\\[-4pt]
{\small University of Oregon}\\[-4pt]
{\small jim@newton.uoregon.edu}
\and
Martin Jackson\\
{\small Department of Mathematics and Computer Science}\\[-4pt]
{\small University of Puget Sound}\\[-4pt]
{\small martinj@ups.edu}
\and
Peng Lu\\
{\small Department of Mathematics}\\[-4pt]
{\small University of Oregon}\\[-4pt]
{\small penglu@darkwing.uoregon.edu}
}
\maketitle

\abstract{ We discuss the Ricci flow on homogeneous 4-manifolds.
After classifying these manifolds, we note that there are families of
initial metrics such that we can diagonalize them and the Ricci flow
preserves the diagonalization.  We analyze the long time behavior of
these families.  We find that if a solution exists for all time, then
the flow exhibits a type III singularity in the sense of Hamilton.  }

\section*{0. Introduction}

It is well-known that there are eight maximal, simply connected
geometries $(X,G)$ with compact quotient in dimension three (\cite{S},
p.474).  In Thurston's geometrization conjecture any closed
three-manifold can be cut into pieces each of which admits one of
these geometries.  To explore the relation between the Ricci flow and
the model geometries, the first two named authors analyze the long
time behavior of the Ricci flow on locally homogeneous three-manifolds
in \cite{IJ}.  In later work (\cite{KM}), using the notion of
quasi-convergence, Knopf and McLeod identify the equivalence classes
of all such flows except the case $X=\widehat{SL}(2,\mathbb{R})$.

Ricci flow has proven to be very successful in studying the geometric
and topological properties of three manifolds (\cite{Ha95}, \cite{P1},
\cite{P2}), and there are indications (\cite{Ha86}, \cite{Ha97},
\cite{Hu2}) that it could be useful for the study of such properties
in four dimensions.  In order to further explore its possible use in
dimension 4 we study the Ricci flow on locally homogeneous
four-manifolds in this paper.  We find that unlike in the case of
three dimensions (\cite{Mi}, \cite{IJ}), some of the families of
locally homogeneous metrics can not be diagonalized because even if
one diagonalizes the initial metric, the flow destroys the
diagonalization of the metric at later times.  In this paper we
identify some families of initial metrics such that the Ricci flow
preserves their diagonalization.  For these families we find that the
behavior of the flow is very close to that seen in dimension three
(\cite{IJ}): either (a) the volume-normalized Ricci flow converges to
a metric of constant sectional curvature or constant holomorphic
bisectional curvature ($\mathbb{C}P^2$ and $\mathbb{C}H^2$); or (b) as
$t \rightarrow +\infty$ the Ricci flow collapses to a lower
dimensional flat manifold with the curvature decaying at the rate
$\frac{1}{t}$; or (c) the Ricci flow approaches, either in finite time
or in infinite time, a direct product of lower dimensional geometries
with constant sectional curvature.

After describing locally homogeneous geometric structures in dimension
4 in section 1, we consider in section 2 the case that the homogeneous
space $X$ is a Lie group .  We identify families of initial metrics
whose diagonalization is preserved by the Ricci flow, and then we
discuss the long time behavior of the Ricci flow for those families.
In section 3, we discuss the long time behavior of the Ricci flow for
the remaining cases.  Since the Ricci flow on closed manifolds
preserves the isometry group, for any locally homogeneous closed
4-manifolds, we discuss the Ricci flow on their universal covering
spaces.

\section*{1. Compact locally homogeneous 4-geometries}

We identify a class of four dimensional homogeneous geometries by
specifying a simply connected four manifold $M$,  a Lie group $G$
that acts transitively on $M$, and the minimal isotropy group $I$ of
the action.  We only consider those $(M,G,I)$ in which $M$ is the
universal cover of a closed manifold $M_q$.  Such a class we call a
\textbf{compact four dimensional homogeneous geometry}.  For each
$(M,G,I)$, there is a collection of Riemannian metrics on $M$ for
which $G$ is the isometry group.  These are the lifts of the locally
homogeneous metrics on $M_q$.  

\subsection*{1.1 List of compact four dimensional homogeneous
geometries}

Let $H^{n}$ be the simply-connected hyperbolic n-manifold and $S^{n}$
be the simply-connected round n-sphere.  We denote the group of
isometries of $H^n$ by $H(n)$.  We summarize the compact three
dimensional homogeneous geometries in the following table.

\medskip

\begin{tabular}{l@{\hspace{24pt}}r@{\hspace{24pt}}r}
Manifold $M^3$ & Lie group $G$ & Isotropy group $I$  \\[2pt]
$\mathbb{R}^3$ & $\mathbb{R}^3$ & $\{0\}$   \\[2pt]
$S^3$ & $SU(2)$ &  $\{ e \}$ \\[2pt]
$\widehat{SL}(2,\mathbb{R})$ & $\widehat{SL}(2,\mathbb{R})$  & $\{ e
\}$\\[2pt]
$Nil^3$ & $Nil^3$  & $\{ e \}$ \\[2pt]
$\widehat{Sol^3}$ & $\widehat{Sol^3}$ & $\{ e\}$ \\[2pt]
$\mathbb{R}^3$ & $E(2)$  &$\{e \}$ \\[2pt]
$S^2 \times \mathbb{R}$ & $SO(3) \times \mathbb{R}$ &  $SO(2) \times
\{0 \}$ \\[2pt]
$H^2 \times \mathbb{R}$ & $H(2) \times \mathbb{R}$  & $SO(2) \times
\{0 \}$ \\[2pt]
$H^3$ & $H(3)$  & $SO(3)$
\end{tabular}

\medskip

\noindent Here $\widehat{SL}(2,\mathbb{R})$ is the universal cover of
the special linear group $SL(2,\mathbb{R})$; its lie algebra $sl_2$
has a basis $X_1,X_2,X_3$ such that the Lie bracket is given by
\[
[X_1,X_2]=-X_3, ~~~ [X_2,X_3]=X_1, ~~~ [X_3,X_1]=X_2.
\]
$Nil^3$ is the 3-dimensional Heisenberg group consisting of
matrices of the form
\[
\left[ 
\begin{array}{rrr}
1 & c_{1} & c_{2} \\ 
0 & 1 & c_{3} \\ 
0 & 0 & 1
\end{array}
\right];
\]
its Lie algebra $n_3$ has a basis $X_1,X_2,X_3$ 
such that the Lie bracket is given by 
\[
[X_1,X_2]=X_3, ~~~ [X_2,X_3]=0, ~~~ [X_3,X_1]=0.
\]
$\widehat{Sol^3}$ is the simply-connected solvable Lie
group 
whose Lie algebra $sol_3$ has a basis $X_1,X_2,X_3$ satisfying
\[
[X_1,X_2]=0, ~~~ [X_2,X_3]= -X_2, ~~~ [X_3,X_1]= -X_1.
\]
$E(2)$ is also a solvable Lie group whose Lie algebra
$L(E_2)$ has a basis $X_1,X_2,X_3$ satisfying
\[
[X_1,X_2]=0, ~~~ [X_2,X_3]=-X_1, ~~~ [X_3,X_1]=-X_2.
\]
 The Lie algebra $su(2)$ of $SU(2)$ can be described by   
\[
[X_1,X_2]=X_3, ~~~ [X_2,X_3]=X_1, ~~~ [X_3,X_1]=X_2.
\]

The compact four dimensional homogeneous geometries have been
classified by Ishihara \cite{I}.  We list them in the following table
(see \cite{W}).


\medskip

\begin{tabular}{r@{\hspace{24pt}}r@{\hspace{24pt}}r}
Manifold $M^4$ & Lie group $G$ & Isotropy group $I$  \\[2pt]
$Nil^4$  &  $Nil^4$   & $\{e\}$ \\[2pt]
$Sol^4_{m,n}$  &  $Sol^4_{m,n}$   & $\{e\}$ \\[2pt]
$Sol_1^4$  &  $Sol_1^4$   & $\{e\}$ \\[2pt]
$Sol_0^4$  &  $Sol_0^4$   & $\{e\}$ \\[2pt]
$\widehat{SL}(2,\mathbb{R}) \times \mathbb{R}$  &
$\widehat{SL}(2,\mathbb{R}) \times \mathbb{R}$   & $\{e\}$ 
\\[2pt]
$Nil^3\times \mathbb{R}$  &  $Nil^3\times \mathbb{R}$   &
$\{e\}$\\[2pt] 
$S^3\times \mathbb{R}$  &  $SU(2)\times\mathbb{R}$   & $\{e\}$ \\[2pt]
$\mathbb{R}^4$  &  $E(2)\times\mathbb{R}$   & $\{e\}$\\[2pt]
$\mathbb{R}^4$  &  $\mathbb{R}^4$   & $\{0\}$\\[2pt]
$S^2\times S^2$  &  $SO(3)\times SO(3)$   & $SO(2)\times SO(2)$\\[2pt]
$S^2\times\mathbb{R}^2$  &  $SO(3)\times \mathbb{R}^2$   &
$SO(2)\times\{0\}$\\[2pt]
$S^2\times H^2$  &  $SO(3)\times H(2)$   & $SO(2)\times SO(2)$\\[2pt]
$ H^2\times\mathbb{R}^2$  &  $H(2)\times\mathbb{R}^2 $   &
$SO(2)\times\{0\} $\\[2pt]
$H^2\times H^2$  &  $H(2)\times H(2)$   & $SO(2)\times SO(2)$\\[2pt]
$\mathbb{C}P^2$  &  $SU(3)$   & $U(2)$\\[2pt]
$\mathbb{C}H^2$  &  $SU(1,2)$   & $U(2)$\\[2pt]
$H^3\times\mathbb{R}$  &  $H(3)\times\mathbb{R}$   &
$SO(3)\times\{0\}$\\[2pt]
$S^4$  &  $SO(5)$   & $SO(4)$\\[2pt]
$H^4$  &  $H(4)$   & $SO(4)$\\[2pt]
\end{tabular}

\medskip

\noindent $Nil^4$, $Sol_{m,n}^4$, $ Sol_1^4$ and $Sol_0^4$ are simply
connected 4-dimensional Lie groups; we describe their Lie algebras in
\S 1.2.  Note that $Sol_{m,n}^4$ includes $\widehat{Sol^3} \times
\mathbb{R}$.  $\mathbb{C}H^2$ is complex hyperbolic space which has
K\"{a}hler symmetric space structure (see \cite{KN}, pp.282-285).

Note that there is another locally homogeneous space $M=F^4$ listed in
\cite{W}.  This is not a compact homogeneous geometry because it does
not have compact quotients.  The isometry group $G$ contains a
discrete subgroup $\Gamma$ such that $F^4/\Gamma$ has finite volume.
  
One can find more detailed description of four dimensional locally
homogeneous geometries in Part II of \cite{Hi}.
  

The Ricci flow study for those classes with trivial isotropy group
requires substantial new analysis; we group these in a category
labelled A.  We describe these classes in \S 1.2.  Those classes with
nontrivial isotropy group are grouped in category B (\S 1.3).

\subsection*{1.2 Four dimensional unimodular Lie groups}

Recall that a Lie group $G$ is called \textbf{co-compact} if $G$
contains a discrete subgroup $\Gamma$ such that $G/\Gamma$ is compact.
Each Lie group in (A) is co-compact.  A co-compact Lie group has
unimodular Lie algebra (\cite{Mi} Lemma 6.2).  Instead of studying
Ricci flow on spaces in (A), we broaden the discussion to Ricci flow
on 4-dimensional unimodular Lie groups.

According to the classification of the 4-dimensional unimodular Lie
algebras (\cite{M}), for each such algebra there is some basis
$X_{1},X_{2},X_{3},X_{4}$ such that the Lie bracket takes the form
indicated below.  We adopt the notation in \cite{M}.

\textbf{A1}. Class $U1[(1,1,1)]$. 
\begin{align*}
\lbrack X_{2},X_{3}] & =0, &  [X_{3},X_{1}]& =0, & \lbrack
X_{1},X_{2}] &  =0,\\
\lbrack X_{1},X_{4}] & =0, &  [X_{2},X_{4}] & =0, & \lbrack
X_{3},X_{4}] &  =0.
\end{align*}
This corresponds to $(M,G,I)=(\mathbb{R}^4,
\mathbb{R}^4,\{0\})$ where $G$ acts on $M$ by translation.

\vskip .3cm

\textbf{A2}. Class $U1[1,1,1]$.
\begin{align*}
\lbrack X_{2},X_{3}] &  =0, & \lbrack X_{3},X_{1}] &  =0, & \lbrack
X_{1},X_{2}] &  =0,\\
\lbrack X_{1},X_{4}] &  =X_{1}, & \lbrack X_{2},X_{4}] &  =kX_{2}, &
\lbrack
X_{3},X_{4}] &  =-(k+1)X_{3},
\end{align*}
where, without loss of generality, we assume $k\geq -\frac{1}{2}$
since otherwise we can interchange $X_2$ and $X_3$.  Only the
following special cases correspond to compact homogeneous geometries.

\noindent \textbf{(A2i)} if $k=0$, the Lie algebra is isomorphic to
the direct sum $sol_3 \oplus \mathbb{R}$.  This corresponds to
$(M,G,I)=(\widehat{Sol^3} \times \mathbb{R}, \widehat{Sol^3} \times
\mathbb{R},\{e\}).$

\noindent \textbf{(A2ii)} if $k=1$, the corresponding geometry is
$(M,G,I)=(Sol_0^4, Sol_0^4,\{e\})$; this can be seen by choosing
$e_1=X_1, e_2=X_2,e_3=X_3$, and $e_4=-X_4$ on p.273 in \cite{W}.

\noindent \textbf{(A2iii)} if there is a number $\alpha >0$ such that
the exponentials of $\alpha$, $\beta \doteq k\alpha$ and $\gamma
\doteq -(k+1)\alpha$ are roots of $\lambda^3 -m\lambda^2 +n \lambda
-1=0$ for some $m,n \in \mathbb{N}$ and $m \neq n$, then one has
$(M,G,I)=(Sol_{m,n}^4, Sol^4_{m,n},\{e\})$ for the geometry.  This can
be seen by choosing $e_1=\alpha X_1, e_2=\alpha X_2,e_3=\alpha X_3$,
and $e_4=-\alpha X_4$ on p.274 and p.270 in \cite{W}.

\vskip .3cm
\textbf{A3}. Class $U1[Z,\bar{Z},1]$.
\begin{align*}
\lbrack X_{2},X_{3}] &  =0, & \lbrack X_{3},X_{1}] &  =0, & \lbrack
X_{1},X_{2}] &  =0,\\
\lbrack X_{1},X_{4}] &  =kX_{1}+X_{2}, & \lbrack X_{2},X_{4}] &
=-X_{1}+kX_{2}, & \lbrack X_{3},X_{4}] &  =-2kX_{3},
\end{align*}
where $k$ is a real number.  if $k=0$, this corresponds to the
geometry $(M,G,I)=(\mathbb{R}^4 , E(2) \times \mathbb{R},\{e\}).$
Other values of $k$ do not correspond to compact homogeneous
geometries.

\vskip .3cm
\textbf{A4}. Class $U1[2,1]$ with $\mu=0$.
\begin{align*}
\lbrack X_{2},X_{3}] &  =0, & \lbrack X_{3},X_{1}] &  =0, & \lbrack
X_{1},X_{2}] &  =0,\\
\lbrack X_{1},X_{4}] &  =X_{2}, & \lbrack X_{2},X_{4}] &  =0, &
\lbrack
X_{3},X_{4}] &  =0.
\end{align*}
This Lie algebra is isomorphic to the direct sum $n_{3}
\oplus\mathbb{R}$ where $n_{3}$ is the Lie algebra of $Nil^3$.  Hence,
in this case $ (M,G,I)=(Nil^3 \times \mathbb{R} , Nil^3 \times
\mathbb{R},\{e\}).$

\vskip .3cm
\textbf{A5}. Class $U1[2,1]$ with $\mu=1$.
\begin{align*}
\lbrack X_{2},X_{3}] &  =0, & \lbrack X_{3},X_{1}] &  =0, & \lbrack
X_{1},X_{2}] &  =0,\\
\lbrack X_{1},X_{4}] &  =-\frac{1}{2}X_{1}+X_{2}, & \lbrack
X_{2},X_{4}] &
=-\frac{1}{2}X_{2}, & \lbrack X_{3},X_{4}] &  =X_{3}.
\end{align*}
This does not correspond to any of the compact homogeneous geometries.

\vskip .3cm
\textbf{A6}. Class $U1[3]$.
\begin{align*}
\lbrack X_{2},X_{3}] &  =0, & \lbrack X_{3},X_{1}] &  =0, & \lbrack
X_{1},X_{2}] &  =0,\\
\lbrack X_{1},X_{4}] &  =X_{2}, & \lbrack X_{2},X_{4}] &  =X_{3}, &
\lbrack
X_{3},X_{4}] &  =0.
\end{align*}
This corresponds to the geometry $(M,G,I)=(Nil^4 , Nil^4,\{e\})$ which
can be seen by choosing $e_1=X_1, e_2=X_2,e_3=X_3$, and $e_4=-X_4$ on
p.274 in \cite{W}.

\vskip .3cm
\textbf{A7}. Class $U3I0$.
\begin{align*}
\lbrack X_{1},X_{4}] &  =0, & \lbrack X_{2},X_{4}] &  =0, & \lbrack
X_{3},X_{4}] &  =0,\\
\lbrack X_{2},X_{3}] &  =X_{4}, & \lbrack X_{3},X_{1}] &  =X_{2}, &
\lbrack
X_{1},X_{2}] &  =-X_{3}.
\end{align*}
This corresponds to the geometry $(M,G,I)=(Sol_1^4 , Sol_1^4,\{e\})$
which can be seen by choosing $e_1=X_1, e_2=X_2+X_3,e_3=X_2 - X_3,e_4=
-2X_4$ on p.272 in \cite{W}.

\vskip .3cm
\textbf{A8}. Class $U3I2$.
\begin{align*}
\lbrack X_{1},X_{4}] &  =0, & \lbrack X_{2},X_{4}] &  =0, & \lbrack
X_{3},X_{4}] &  =0,\\
\lbrack X_{2},X_{3}] &  =-X_{4}, & \lbrack X_{3},X_{1}] &  =X_{2}, &
\lbrack
X_{1},X_{2}] &  =X_{3}.
\end{align*}
This does not correspond to any of the compact homogeneous geometries.

\vskip .3cm
\textbf{A9}. Class $U3S1$.
\begin{align*}
\lbrack X_{1},X_{4}] &  =0, & \lbrack X_{2},X_{4}] &  =0, & \lbrack
X_{3},X_{4}] &  =0,\\
\lbrack X_{2},X_{3}] &  =X_{1}, & \lbrack X_{3},X_{1}] &  =X_{2}, &
\lbrack
X_{1},X_{2}] &  =-X_{3}.
\end{align*}
This Lie algebra is isomorphic to the direct sum $sl_{2}
\oplus\mathbb{R}$.  This corresponds to the geometry
$(M,G,I)=(\widehat{SL}(2,\mathbb{R}) \times \mathbb{R},
\widehat{SL}(2,\mathbb{R}) \times \mathbb{R},\{e\}).$

\vskip .3cm
\textbf{A10}. Class $U3S3$.
\begin{align*}
\lbrack X_{1},X_{4}] &  =0, & \lbrack X_{2},X_{4}] &  =0, & \lbrack
X_{3},X_{4}] &  =0,\\
\lbrack X_{2},X_{3}] &  =X_{1}, & \lbrack X_{3},X_{1}] &  =X_{2}, &
\lbrack
X_{1},X_{2}] &  =X_{3}.
\end{align*}
This Lie algebra is isomorphic to the direct sum $su(2)\oplus
\mathbb{R}$.  This corresponds to the geometry $(M,G,I)=(S^3 \times
\mathbb{R}, SU(2) \times \mathbb{R},\{e\}).$

\subsection*{1.3 Compact four dimensional homogeneous geometries 
with nontrivial isotropy group}

Now we list the compact 4-dimensional homogeneous geometries
$(M^4,G,I)$ for which dimension of $G$ is bigger than 4.  Recall
$H(n)$ is the isometry group of the simply-connected hyperbolic
n-manifolds $H^n$.

\vskip .2cm
B1. $(M,G,I) = (H^{3} \times \mathbb{R},H(3)\times \mathbb{R}, SO(3))$

\vskip .2cm
B2. $(M,G,I) = (S^2 \times \mathbb{R}^{2}, SO(3) \times \mathbb{R}^2, 
SO(2) \times \{0 \})$

\vskip .2cm
B3. $(M,G,I) = (H^2 \times \mathbb{R}^2, H(2) \times \mathbb{R}^2,
SO(2) 
\times \{0 \}) $

\vskip .2cm
B4. $(M,G,I) =(S^2 \times S^2, SO(3) \times SO(3), SO(2) \times 
SO(2))$
 
\vskip .2cm
B5. $(M,G,I) = (S^{2} \times H^{2},  SO(3)\times H(2), 
SO(2) \times SO(2))$

\vskip .2cm
B6. $(M,G,I) = H^2 \times H^2, H(2) \times H(2),
SO(2) \times SO(2))$

\vskip .2cm
B7. $(M,G,I) = (\mathbb{C}P^2, SU(3), U(2))$

\vskip .2cm
B8. $(M,G,I) = (\mathbb{C}H^2, SU(1,2), U(2))$ 

\vskip .2cm
B9. $(M,G,I) = (S^4, SO(5), SO(4))$ 

\vskip .2cm
B10. $(M,G,I) =(H^4, H(4), SO(4))$

\section*{2. The Ricci flow on 4-dimensional unimodular Lie groups}

Recall our strategy is to analyze the Ricci flow on a simply connected
manifold $M$ that is the universal cover of a closed manifold $M_q$.
For a fixed class $(M,G,I)$ and a fixed initial homogeneous metric
$g_0$ compatible with the class, let $g(t)$ be the homogeneous
solution of the Ricci flow
\[
\frac{\partial}{\partial t}
{g}(t)=-2\operatorname{Ric}({g}(t)) \qquad g(0)=g_0 .
\]
On the closed manifold $M_q$, we consider the volume-normalized Ricci
flow ${g}_N(t)$ as the solution of
\[
\frac{\partial}{\partial t}
{g}_N(t)=-2\operatorname{Ric}({g}_N(t))+\frac{r_N}{2}{g}_N\qquad
g_N(0)=g_0
\]
where $r_N$ is the scalar curvature of $g_N$.  Note that averaging of
the scalar curvature is not needed for homogeneous metrics.  We also
note that the Ricci flow equation for homogeneous metrics reduces to a
system of ordinary differential equations.

For each class (Ai) listed in \S 1.2, we describe the families of
initial metrics which are diagonal and remain diagonal under the Ricci
flow and then study their long time behavior.  To address the
diagonalization issue, we use the following strategy: Fix a
homogeneous metric $h$ and let $\{X_{i}\}$ be any basis of
left-invariant vectors on the Lie group $G$ with the bracket structure
\[
[X_{i},X_{j}]=C^{k}_{\phantom{k} ij}X_{k}.
\]
For those classes (e.g., A4, A9 and A10) in which the Lie group 
$G$ is a product group $G_1 \times \mathbb{R}$ with $\dim(G_1)=3$, we
choose $X_1,X_2,X_3$ as left invariant vector fields on $G_1$ and 
$X_4=\frac{\partial}{\partial u}$ on $\mathbb{R}$.

Let $\{Y_{i}\}$ be any basis orthogonal with respect to $h$; 
that is,
\[
h(Y_{i},Y_{j})=\lambda_{i}\delta_{ij}.
\]
Let the transformation from $\{X_{i}\}$ to $\{Y_{i}\}$ be given by
\[
Y_{i}=\Lambda^{k}_{\phantom{k} i}X_{k}.
\]
Computing the bracket structure for $\{Y_{i}\}$, we get
\[
[Y_{i},Y_{j}]=[\Lambda^{k}_{\phantom{k}i}X_{k},
\Lambda^{l}_{\phantom{k}j}X_{l}]
=\Lambda^{k}_{\phantom{k} i}\Lambda^{l}_{\phantom{k}j}
C^{m}_{\phantom{k}kl}X_{m} 
=\Lambda^{k}_{\phantom{k} i}\Lambda^{l}_{\phantom{k} j}
C^{m}_{\phantom{k} kl} (\Lambda^{-1})^{n}_{\phantom{k} m}Y_{n}.
\]
Thus,
\[
[Y_{i},Y_{j}]=\tilde{C}^{n}_{\phantom{k} ij}Y_{n}
\]
where
\[
\tilde{C}^{n}_{\phantom{k} ij}= \Lambda^{k}_{\phantom{k} i}\Lambda
^{l}_{\phantom{k} j}(\Lambda^{-1})^{n}_{\phantom{k}
m}C^{m}_{\phantom{k} kl}.
\]

Now compute the Ricci curvature of $h$ using the orthonormal basis
$\{\bar{Y}_{i}\}$ defined by
$\bar{Y}_{i}=\frac{1}{\sqrt{\lambda_{i}}}Y_{i}$.  For this, we use the
following Ricci curvature formula for unimodular Lie groups from
Corollary 7.33 p.184 \cite{B},

\begin{equation*}
\operatorname{Ric}(W,W)=
-\frac{1}{2}\sum_{i}\bigl|[W,\bar{Y}_{i}]\bigr|^{2}
-\frac{1}{2}\sum_{i}
\langle\bigl[W,[W,\bar{Y}_{i}]\bigr],\bar{Y}_{i}\rangle\\
+\frac{1}{2}\sum_{i<j}\langle[\bar{Y}_{i},\bar{Y}_{j}],W\rangle^{2} .
 \tag{1}  \label{eq 1}
\end{equation*}
Finally, check if any positive values of the parameters $\lambda_{i}$
produce a diagonal Ricci tensor.  Only for these values does the
metric remain diagonal under the Ricci flow.  We follow this strategy
and use the same notation in the rest of this section.

\textbf{Remark}.  In our search for families of initial metrics which
remain diagonal under the Ricci flow, we have chosen special $Y_i$ so
that the Lie brackets $[Y_i,Y_j]$ are simple.  Presumably there are
other families of initial metrics and other bases $Y_i$ for which the
diagonalization is preserved by the Ricci flow.  Our calculations show
there are $(M,G,I)$ and bases $Y_i$ such that the property that the
initial metric has components $(g_0)_{a4}=0$, $a=1,2,3$, is preserved
under Ricci flow.

To study the decay of the curvature tensor, we use the following
sectional curvature formula for Lie groups from Theorem 7.30 p.183
\cite{B}.  For the Lie algebra $\mathfrak{g}$ of $G$, define the
operator $U:\mathfrak{g} \times \mathfrak{g} \rightarrow \mathfrak{g}$
by
\[
2\langle U(X,Y),Z\rangle=\langle[Z,X],Y\rangle+\langle X,[Z,Y]\rangle 
\text{ for all } Z \in \mathfrak{g};
\tag{2} \label{eq 2} 
\]
then the curvature is given by
\begin{multline*}
\langle R(X,Y)X,Y\rangle=
-\frac{3}{4} | [X,Y]|^2 
-\frac{1}{2}\langle[X,[X,Y]],Y\rangle 
-\frac{1}{2}\langle[Y,[Y,X]],X\rangle \\
+|U(X,Y) |^2 
-\langle U(X,X),U(Y,Y)\rangle. \tag{3} \label{eq 3} 
\end{multline*}

\subsection*{A1. U1[(1,1,1)]}

For $(M,G,I)=(\mathbb{R}^{4},\mathbb{R}^{4},\{0\})$, $G$ acts on $M$
by translation $h(x)=h+x$ for $h\in G$.  Any homogeneous metric $g_0$
on $M$ must be of the form
\[
g_0 = \lambda_1 dx^1 \otimes dx^1+
+ \lambda_2 dx^2 \otimes dx^2
+ \lambda_3 dx^3 \otimes dx^3
+ \lambda_4 dx^4 \otimes dx^4
\]
for some constants $ \lambda_i>0$.  The metric $g_0$ is flat; hence
$g(t)\equiv g_{0}$ for $-\infty<t<\infty$.

\subsection*{A2. $U1[1,1,1]$}

For $U1[1,1,1]$, we use $Y_{i}=\Lambda^{k}_{\phantom{k} i}X_{k}$ with
\[
\Lambda=
\begin{bmatrix}
1 & 0 & 0 & 0\\
a_{1} & 1 & 0 & 0\\
a_{2} & a_{3} & 1 & 0\\
a_{4} & a_{5} & a_{6} & 1
\end{bmatrix}
\]
to diagonalize the initial metric $g_0$.

\begin{proposition}
For the class $U1[1,1,1]$ suppose the initial metric $g_{0}$ is
diagonal in the basis $Y_i$. Then

(i) if $k\neq 1$ the Ricci flow solution $g(t)$ remains diagonal 
in the basis $Y_i$ if and only if $a_1=a_2=a_3=0$;

(ii) if $k=1$ the Ricci flow solution $g(t)$ remains diagonal 
in the basis $Y_i$ if and only if $a_2=a_3=0$.
\end{proposition}

\begin{proof} We compute
\begin{align*}
  [Y_{2},Y_{3}]  &  =0, 
& [Y_{3},Y_{1}]  &  =0, 
& [Y_{1},Y_{2}]  &  =0,\\
  [Y_{1},Y_{4}]  &  =Y_{1}, 
& [Y_{2},Y_{4}]  &  =kY_{2}+\alpha Y_{1}, 
& [Y_{3},Y_{4}]  &  =-(k+1)Y_{3}+\beta Y_{2}+\gamma Y_{1}.
\end{align*}
where
\[
\alpha=(1-k)a_{1},\quad\beta=(1+2k)a_{3},\quad\text{and}\quad\gamma
=(k+2)a_{2}-(1+2k)a_{1}a_{3}.
\]

Let $W=w_{1} \bar{Y}_{1} +w_{2} \bar{Y}_{2} + w_{3} \bar{Y}_{3} +w_{4}
\bar{Y}_{4}$.  We compute $[W,\bar{Y}_i]$ first and then use (\ref{eq
1}) with $h=g_0$ to compute the coefficients of $w_iw_j$ in
$\operatorname{Ric}(W,W)$.  We find that the off-diagonal components
of the Ricci tensor in the basis $\{\bar{Y}_{i}\}$ are given by
\begin{align*}
& \operatorname{Ric}(\bar{Y}_{1},\bar{Y}_{2})= \frac{\bigl(\beta
\gamma\lambda_{2}+(k-1)\alpha\lambda_{3}\bigr)\sqrt{\lambda_{1}}}
{2\sqrt{\lambda_{2}}\lambda_{3}\lambda_{4}}\\
& \operatorname{Ric}(\bar{Y}_{1},\bar{Y}_{3})= -\frac{(2+k)\gamma
\sqrt{\lambda_{1}}} {2\sqrt{\lambda_{3}}\lambda_{4}}\\
&\operatorname{Ric}(\bar{Y}_{2},\bar{Y}_{3})= -\frac{\alpha
\gamma\lambda_{1}+(1+2k)\beta\lambda_{2}} {2\sqrt{\lambda_{2}}\sqrt
{\lambda_{3}}\lambda_{4}}\\
&  \operatorname{Ric}(\bar{Y}_{1},\bar{Y}_{4})=
\operatorname{Ric}(\bar{Y}_{2},\bar{Y}_{4})=
\operatorname{Ric}(\bar{Y}_{3},\bar{Y}_{4})=0
\end{align*}

Recall $k\geq -\frac{1}{2}$.  Note that we have $\beta=0$ by
definition for $k=-\frac{1}{2}$.  In order for these off-diagonal
components to be zero, we must have $\alpha=\beta=\gamma=0$ if $k\neq
1$ and $\beta=\gamma=0$ if $k= 1$.
\end{proof}

Now we discuss the long time behavior for the families of the initial
metrics in Proposition 1.  We start with (A2iii) and later show that
the other two cases are covered by the same analysis.

\textbf{(A2iii)}. In this case $k \neq 0,1$. We have 
\begin{align*}
& Y_1=X_1, & Y_2=& X_2, \\
& Y_3=X_3, & Y_4=& X_4+a_4X_1+a_5X_2+a_6X_3.
\end{align*}
and 
\begin{align*}
[Y_{2},Y_{3}]  &  =0, & [Y_{3},Y_{1}]  &  =0, & [Y_{1},Y_{2}]  &
=0,\\
[Y_{1},Y_{4}]  &  =Y_{1}, & [Y_{2},Y_{4}]  &  =kY_{2}, &
[Y_{3},Y_{4}]  &  =-(k+1)Y_{3}.
\end{align*}

The bases $Y_i$ and $X_i$ both satisfy the same Lie bracket relations
so either can be used in the Ricci flow analysis.  We use $X_i$.

Let $\theta_{i}$ be the frame of 1-forms dual to $X_{i}$. 
Assume the Ricci flow solution takes the special form
\[
 g(t)= A(t) (\theta_{1})^{2} +B(t) (\theta_{2})^{2} + 
 C(t)(\theta_{3})^{2} +D(t) (\theta_{4})^{2}
\]
with
\[
g_0 = \lambda_1 (\theta_{1})^{2} + \lambda_2 (\theta_{2})^{2} + 
\lambda_3 (\theta_{3})^{2} + \lambda_4 (\theta_{4})^{2}.
\]
Then
$\bar{X}_1=\frac{1}{\sqrt{A}}X_1,\cdots,\bar{X}_4=\frac{1}{\sqrt{D}}X_4$
is an orthonormal frame with respect to the metric $g$.
Let $W=w_{1} \bar{X}_{1} +w_{2} \bar{X}_{2} + w_{3} \bar{X}_{3} 
+w_{4} \bar{X}_{4}$ and then compute
\begin{align*}
&  [W,\bar{X}_{1}]= - \frac{1}{\sqrt{D}} w_{4} \bar{X}_{1}
&  [W,\bar{X}_{2}]= &- \frac{k}{\sqrt{D}} w_{4} \bar{X}_{2}\\
&  [W,\bar{X}_{3}]= \frac{k+1}{\sqrt{D}} w_{4} \bar{X}_{3}
&  [W,\bar{X}_{4}]= &\frac{1}{\sqrt{D}} w_{1} \bar{X}_{1} +
\frac{k}{\sqrt{D}}
w_{2} \bar{X}_{2} - \frac{k+1}{\sqrt{D}} w_{3} \bar{X}_{3}.
\end{align*}
We have from (\ref{eq 1}) with $h=g$
\begin{align*}
\operatorname{Ric}(W,W)= 0 \cdot w_{1}^{2} +0 \cdot w_{2}^{2}+ 0
\cdot w_{3}^{2} -\frac{2(k^{2}+k+1)}{D} \cdot w_{4}^{2}.
\end{align*}
So
\begin{align*}
&  \operatorname{Ric}(X_{1},X_{1} )= \operatorname{Ric}(X_{2},X_{2} )=
\operatorname{Ric}(X_{3},X_{3} )= 0,\\
&  \operatorname{Ric}(X_{4},X_{4} )= D
\cdot\operatorname{Ric}(\bar{X}%
_{4},\bar{X}_{4} )= -2(k^{2}+k+1),
\end{align*}
and the Ricci flow is
\begin{align*}
&  \frac{dA}{dt}=0,
&  \frac{dB}{dt}=&0,\\
&  \frac{dC}{dt}=0,
&  \frac{dD}{dt}=& 4 (k^{2}+k+1).
\end{align*}
The solution is given by
\[
A(t)=\lambda_{1},\quad B(t)=\lambda_{2},\quad
C(t)=\lambda_{3},\quad\text{and}\quad D(t)=\lambda_{4}+4
(k^{2}+k+1)t.
\]

Hence for the subfamily in Proposition 1, the Ricci flow does not move
in three directions and expands in the fourth direction at a speed
linear in $t$.  Pick a point $p \in M_q$.  It is clear that the
volume-normalized solution $(M_q,{g}_N(t),p)$ converges/collapses to a
line in the the pointed Gromov-Hausdorff topology.
  
Next we compute the curvature decay of $g(t)$. From (\ref{eq 2}) we
find
\begin{xalignat*}{3}
   U(X_1,X_1) & = -\frac{A}{D}X_4 
&  U(X_2,X_2) & =  -\frac{kB}{D}X_4
&  U(X_3,X_3) & =  \frac{(k+1)C}{D}X_4   \\
 U(X_4,X_4) & =  0
& U(X_1,X_2) & =  0 
& U(X_1,X_3) & =  0 \\
 U(X_2,X_3) & =  0 
& U(X_1,X_4) & =  \frac{1}{2}X_1 
& U(X_2,X_4) & = \frac{k}{2}X_2 \\
 U(X_3,X_4) & =  -\frac{k+1}{2}X_3. 
\end{xalignat*} 
From (\ref{eq 3}) with $h=g$ we find the sectional
curvatures
\begin{xalignat*}{3}
  K(X_1,X_2) & =-\frac{k}{D}, 
& K(X_1,X_3) & =\frac{k+1}{D},
& K(X_2,X_3) & =\frac{k(k+1)}{D},   \\
  K(X_1,X_4) & =-\frac{1}{D}, 
& K(X_2,X_4) & = -\frac{k^2}{D}, 
& K(X_3,X_4) & =-\frac{(k+1)^2}{D}.
\end{xalignat*} 
These curvatures of the solution $g(t)$ decay at the rate $1/t$.

\textbf{(A2i)}.  This is a special case of (A2iii) if we allow $k=0$,
so the analysis in (A2iii) applies.  Note that since the Lie algebra
is the direct sum $sol_3 \oplus \mathbb{R}$, we can get the same
conclusions from the analysis in \cite{IJ}(pp.733-735).

\textbf{(A2ii)}. In this case $k=1$. We have 
\begin{align*}
& Y_1=X_1, & Y_2=& X_2+a_1X_1, \\
& Y_3=X_3, & Y_4=& X_4+a_4X_1+a_5X_2+a_6X_3.
\end{align*}
and 
\begin{align*}
[Y_{2},Y_{3}]  &  =0, & [Y_{3},Y_{1}]  &  =0, & [Y_{1},Y_{2}]  &
=0,\\
[Y_{1},Y_{4}]  &  =Y_{1}, & [Y_{2},Y_{4}]  &  =Y_{2}, &
[Y_{3},Y_{4}]  &  =-2Y_{3}.
\end{align*}
Note that this is the Lie algebra structure in (A2iii) if we allow
$k=1$.  Hence the analysis of (A2iii) applies with the same
conclusion.


\subsection*{A3. $U1[Z,\bar{Z},1]$}

For $U1[Z,\bar{Z},1]$, we use $Y_{i}=\Lambda^{k}_{\phantom{k}
i}X_{k}$ with
\[
\Lambda=
\begin{bmatrix}
1 & a_{2} & a_{3} & 0\\
0 & 1 & a_{1} & 0\\
0 & 0 & 1 & 0\\
a_{4} & a_{5} & a_{6} & 1
\end{bmatrix}
\]
to diagonalize the initial metric $g_0$.

\begin{proposition}
For the class $U1[Z,\bar{Z},1]$ suppose the initial metric
$g_0$ is diagonal in the basis $Y_i$. 
Then the Ricci flow solution $g(t)$ remains diagonal in the basis $Y_i$
if and only if $a_1=a_2=a_3=0$.
\end{proposition}

\begin{proof} We compute
\begin{align*}
&  [Y_{2},Y_{3}]=0, \qquad[Y_{3},Y_{1}]=0,
\qquad[Y_{1},Y_{2}]=0,\qquad 
[Y_{3}, Y_{4}]=-2kY_{3} \\
&  [Y_{1},Y_{4}]=(k-\alpha)Y_{1}+(\alpha^{2}+1)Y_{2}+\beta Y_{3}, 
~~ [Y_{2},Y_{4}]=-Y_{1}+(k+\alpha)Y_{2}+\gamma Y_{3}.
\end{align*}
where
\[
\alpha=a_{2},\quad\beta=a_{2}a_{3}-a_{1}a_{2}^{2}-a_{1}-3ka_{3},\quad
\text{and}\quad\gamma=a_{3}-3ka_{1}-a_{1}a_{2}.
\]

We compute the off-diagonal components of the Ricci tensor in the
basis $\{\bar{Y}_{i}\}$ using (\ref{eq 1}) as in \S 2.A2 and get
\begin{align*}
& \operatorname{Ric}(\bar{Y}_{1},\bar{Y}_{2})= -\frac{2\alpha
\lambda_{1}+2\alpha(1+\alpha^{2})\lambda_{2} +\beta\gamma\lambda_{3}}%
{2\sqrt{\lambda_{1}\lambda_{2}}\lambda_{4}}\\
&  \operatorname{Ric}(\bar{Y}_{1},\bar{Y}_{3})= -\frac{\bigl(\gamma
\lambda_{1}+(\alpha-3k)\beta\lambda_{2}\bigr)\sqrt{\lambda_{3}}}
{2\sqrt{\lambda_{1}}\lambda_{2}\lambda_{4}}\\
&  \operatorname{Ric}(\bar{Y}_{2},\bar{Y}_{3})= \frac{\bigl((\alpha
+3k)\gamma\lambda_{1}+(1+\alpha^{2})\beta\lambda_{2}\bigr)\sqrt{\lambda_{3}}}
{2\lambda_{1}\sqrt{\lambda_{2}}\lambda_{4}}\\
&   \operatorname{Ric}(\bar{Y}_{1},\bar{Y}_{4})=
\operatorname{Ric}(\bar{Y}_{2},\bar{Y}_{4})=
\operatorname{Ric}(\bar{Y}_{3},\bar{Y}_{4})=0.
\end{align*}
In order for these off-diagonal components to be zero, we must have
$\alpha=\beta=\gamma=0$ and the proposition follows.
\end{proof}

If $\alpha=\beta=\gamma=0$, the $Y_i$ and $X_i$ both satisfy the same
Lie bracket relations.  As in \S 2.A2, we use $X_i$ in carrying out
the analysis of the long time behavior of the Ricci flow solution for
the family in Proposition 2.  Proceeding as in \S 2.A2, we find
\begin{align*}
&  [W,\bar{X}_{1}]= - \frac{k}{\sqrt{D}} w_{4} \bar{X}_{1}
-\sqrt{\frac{B}%
{AD}} w_{4} \bar{X}_{2}\\
&  [W,\bar{X}_{2}]= \sqrt{\frac{A}{BD}} w_{4} \bar{X}_{1} - \frac
{k}{\sqrt{D}} w_{4} \bar{X}_{2}, \qquad  [W,\bar{X}_{3}]=
\frac{2k}{\sqrt{D}} w_{4} \bar{X}_{3}\\
&  [W,\bar{X}_{4}]= ( \frac{k}{\sqrt{D}} w_{1} - \sqrt{\frac{A}{BD}} 
w_{2} ) \bar{X}_{1} + (\sqrt{\frac{B}{AD}} w_{1} + \frac{k}{\sqrt{D}}
w_{2} ) \bar{X}_{2} - \frac{2k}{\sqrt{D}} w_{3} \bar{X}_{3}.
\end{align*}
We have from (\ref{eq 1}) 
\begin{align*}
\operatorname{Ric}(W,W) = \frac{A^2-B^2}{2ABD}w_{1}^{2} 
-\frac{A^2-B^2}{2ABD} w_{2}^{2} + 0 \cdot w_{3}^{2} - \frac{
(A-B)^2+12k^{2} AB}{2ABD}w_{4}^{2},
\end{align*}
so the Ricci flow is
\begin{xalignat*}{2}
  \frac{dA}{dt}  & =- \frac{A^{2} -B^{2}}{BD}, 
&  \frac{dB}{dt} & =-\frac{B^{2} -A^{2}}{AD},\\
  \frac{dC}{dt}  & =0, 
&  \frac{dD}{dt} & = \frac{(A-B)^{2}+ 12k^{2} AB}{AB}.
\end{xalignat*}
Clearly $C(t)=\lambda_3$.

If $\lambda_{1}=\lambda_{2}$, then
\[
A(t)=\lambda_1, \quad B(t)=\lambda_2,\quad\text{and} 
\quad D(t)=\lambda_4+12k^{2}t.
\]

If $\lambda_{1} \neq \lambda_{2}$, we may assume $\lambda_{2}>
\lambda_{1}$ without loss of generality by the symmetry of $A$ and $B$
in this system.  A simple computation gives
\[
\frac{1}{A}\frac{dA}{dt}+\frac{1}{B}\frac{dB}{dt}=0,
\]
so the product $AB=\lambda_1 \lambda_2$ for all $t$. Another
computation gives
\[
\frac{d}{dt}[A-B]=-\frac{(A+B)^{2}}{ABD}(A-B),
\]
so $A-B$ is decreasing in $t$.  From the equation for $\frac{dD}{dt}$
it follows easily that
\[
12k^{2} \leq \frac{dD}{dt} \leq 12k^{2}+\frac{\lambda_2}{\lambda_1},
\]
and so
\[
\lambda_{4}+12k^{2}t \leq D(t) \leq \lambda_{4}+
(12k^{2}+\frac{\lambda_2}{\lambda_1})t .
\]

Hence for the family in Proposition 2, the long-time behavior of the
solution $g(t)$ as $t \rightarrow +\infty$ is
\[
A(t)\to\sqrt{\lambda_1 \lambda_2},\quad B(t)\to\sqrt{\lambda_1
\lambda_2},
\quad C(t)=\lambda_3, \quad D(t)\to\infty \text{ linearly}.
\] 
Pick a point $p \in M_q$.  It is clear that the volume-normalized
solution $(M_q,{g}_N(t),p)$ converges/collapses to a line in the pointed
Gromov-Hausdorff topology.
  
Next we compute the curvature decay of $g(t)$. From (\ref{eq 2}) we
find
\begin{xalignat*}{3}
  U(X_1,X_1) & =-\frac{kA}{D}X_4
& U(X_2,X_2) & =-\frac{kB}{D}X_4
& U(X_3,X_3) & =\frac{2kC}{D}X_4   \\
  U(X_4,X_4) & =0 
& U(X_1,X_2) & =\frac{A-B}{2D}X_4 
& U(X_1,X_3) & =0 \\
  U(X_2,X_3) & =0
& U(X_1,X_4) & = \frac{k}{2}X_1-\frac{A}{2B}X_2 
& U(X_2,X_4) & =\frac{B}{2A}X_1+\frac{k}{2}X_2 \\
  U(X_3,X_4) & =-kX_3.  
\end{xalignat*} 
From (\ref{eq 3}) with $h=g$ we find the sectional
curvatures
\begin{xalignat*}{2}
  K(X_1,X_2) & =\frac{\frac{A}{B}+\frac{B}{A}-2-4k^2}{4D} 
& K(X_1,X_3) & =\frac{2k^2}{D}
\\
 K(X_2,X_3) & =\frac{2k^2}{D}
& K(X_1,X_4) & =\frac{\frac{A}{B}-3\frac{B}{A}+2-4k^2}{4D}
\\
K(X_2,X_4) & = \frac{-{3}\frac{A}{B}+\frac{B}{A}+{2}-4k^2}{D} 
& K(X_3,X_4) & =-\frac{4k^2}{D}.
\end{xalignat*} 
Hence for the family in Proposition 2, 
the curvatures of the solution $g(t)$ decay at the rate $1/t$.

\subsection*{A4. $U1[2,1]$}
For $U1[2,1]$, we use $Y_{i}=\Lambda^{k}_{\phantom{k} i}X_{k}$ with
\[
\Lambda=
\begin{bmatrix}
1 & a_{2} & a_{3} & 0\\
0 & 1 & 0 & 0\\
0 & a_{1} & 1 & 0\\
a_{4} & a_{5} & a_{6} & 1
\end{bmatrix}
\]
to diagonalize the initial metric $g_0$.

\begin{proposition}
For the class $U1[2,1]$ suppose the initial metric $g_{0}$ is diagonal
in the basis $Y_i$. Then the Ricci flow solution $g(t)$ remains diagonal
in the basis $Y_i$ for all $t \geq 0$.
\end{proposition}

\begin{proof} We compute
\begin{align*}
[Y_{2},Y_{3}]  &  =0, & [Y_{3},Y_{1}]  &  =0, & [Y_{1},Y_{2}]  &
=0,\\
[Y_{1},Y_{4}]  &  =Y_{2}, & [Y_{2},Y_{4}]  &  =0, & [Y_{3},Y_{4}]  &
=0.
\end{align*}
This bracket structure is identical to that of the basis $\{X_{i}\}$. 

We compute the off-diagonal components of the Ricci tensor in the
basis $\{ \bar{Y}_i\}$ using (\ref{eq 1}) as in \S 2.A2 and find
$\operatorname{Ric}(\bar{Y}_{i},\bar{Y}_{j})=0$ for all $i<j$.  The
proposition is proved.
\end{proof}

Since $Y_i$ and $X_i$ both satisfy the same Lie bracket relations, as
in \S 2.A2 we use $X_i$ and can carry out the analysis of the long
time behavior of the Ricci flow solution for the family in Proposition
3.  Proceeding as in \S 2.A2, we find the Ricci tensor
\[
\operatorname{Ric}(W,W)= -\frac{B}{2AD}w_1^2+ \frac{B}{2AD}w_2^2
+0\cdot w_3^2-\frac{B}{2AD}w_4^2.
\]
Hence the Ricci flow is
\begin{align*}
&  \frac{dA}{dt}= \frac{B}{D}, & & \frac{dB}{dt}=- \frac{B^2}{AD},\\
&  \frac{dC}{dt}=0, & & \frac{dD}{dt}= \frac{B }{A}.
\end{align*}

From $\frac{d}{dt}(\frac{A}{D})=\frac{d}{dt}(AB)=0$ we get
\begin{align*}
&A=(\lambda_1+\frac{3\lambda_1^2\lambda_2}{\lambda_3}t)^{1/3} & 
& B= \lambda_1\lambda_2 {(\lambda_1+\frac{3\lambda_1^2\lambda_2
}{\lambda_3}t)^{-1/3}} \\
& C=\lambda_3 & & D=\frac{\lambda_4}{\lambda_1}(\lambda_1+
\frac{3\lambda_1^2\lambda_2}{\lambda_3}t)^{1/3}
\end{align*}
Hence for the family in Proposition 3, the long time behavior of the
solution $g(t)$ as $t \to +\infty$ is
\[
A(t) \to +\infty, \qquad B(t) \to 0^+, \qquad C(t)=\lambda_3,
\qquad D(t) \to +\infty.
\]

Pick a point $p \in M_q$.  It is clear that the volume-normalized
solution $(M_q,{g}_N(t),p)$ converges/collapses to a plane in the
pointed Gromov-Hausdorff topology.

Next we compute the curvature decay for $g(t)$.  From (\ref{eq 2}) we
find $U(X_1,X_2)=-\frac{B}{2D}X_4$, $U(X_2,X_4)=\frac{B}{2A}X_1$ and
all other $U(X_i,X_j)=0$.  From (\ref{eq 3}) with $h=g$ we find the
sectional curvatures
\[
K(X_1,X_2)=K(X_2,X_4)=\frac{B}{4AD}, \qquad
K(X_1,X_4)=-\frac{3B}{4AD},
\]
and all other $K(X_i,X_j)=0$.  Hence for the family in Proposition 3,
the curvatures of the solution $g(t)$ decay at the rate $1/t$.  Note
that since the Lie algebra is the direct sum $n_3 \oplus \mathbb{R}$,
we can get the same conclusions from the analysis in \cite{IJ}(p.734).

\subsection*{A5. $U1[2,1]$}

For $U1[2,1]$, we use $Y_{i}=\Lambda^{k}_{\phantom{k} i}X_{k}$ with
\[
\Lambda=
\begin{bmatrix}
1 & a_{2} & a_{3} & 0\\
0 & 1 & a_{1} & 0\\
0 & 0 & 1 & 0\\
a_{4} & a_{5} & a_{6} & 1
\end{bmatrix}
\]
to diagonalize the initial metric $g_0$.

\begin{proposition}
For the class $U1[2,1]$ suppose the initial metric $g_{0}$ is diagonal
in the basis $Y_i$.  Then the Ricci flow solution $g(t)$ remains
diagonal in the basis $Y_i$ if and only if $a_1=a_3=0$.
\end{proposition}

\begin{proof} We compute
\begin{align*}
[Y_{2},Y_{3}]  &  =0, & [Y_{3},Y_{1}]  &  =0, & [Y_{1},Y_{2}]  &
=0,\\
[Y_{1},Y_{4}]  &  =-\frac{1}{2}Y_{1}+Y_{2}+\beta Y_{2}, &
[Y_{2},Y_{4}]  &
=-\frac{1}{2}Y_{2}+\alpha Y_{3}, & [Y_{3},Y_{4}]  &  =Y_{3}.
\end{align*}
where
\[
\alpha=\frac{3}{2}a_{1} \qquad \beta=\frac{3}{2}(a_{3}-a_{1}).
\]

We compute the off-diagonal components of the Ricci tensor in the
basis $\{\bar{Y}_{i}\}$ using (\ref{eq 1}) as in \S 2.A2 and get
\begin{xalignat*}{3}
\operatorname{Ric}(\bar{Y}_{1},\bar{Y}_{2})
& =-\frac{\alpha\beta\lambda_{3}}{2\sqrt{\lambda_{1}\lambda_{2}}\lambda_{4}}
& \operatorname{Ric}(\bar{Y}_{1},\bar{Y}_{3})
& = -\frac{3\beta\sqrt{\lambda_{3}}}{4\sqrt{\lambda_{1}}\lambda_{4}}
& \operatorname{Ric}(\bar{Y}_{1},\bar{Y}_{4})
& =0\\
\operatorname{Ric}(\bar{Y}_{2},\bar{Y}_{3})
& = \frac{(-3\alpha\lambda_{1}+2\beta\lambda_{2})\sqrt{\lambda_{3}}} 
{4\lambda_{1}\sqrt{\lambda_{2}}\lambda_{4}} 
& \operatorname{Ric}(\bar{Y}_{2},\bar{Y}_{4}) & =0
& \operatorname{Ric}(\bar{Y}_{3},\bar{Y}_{4}) & =0.
\end{xalignat*}
In order for these off-diagonal components to be zero, we must have
$\alpha=\beta=0$ and the proposition follows.
\end{proof}

If $\alpha=\beta=0$, the $Y_i$ and $X_i$ both satisfy the same Lie
algebra bracket relations.  As in \S 2.A2, we use $X_i$ and carry out
the analysis of the long time behavior of the Ricci flow solution for
the family in Proposition 4.  Proceeding as in \S 2.A2, we find
 \begin{align*}
[W,\bar{X}_{1}] 
& = \frac{1}{2\sqrt{D}} w_{4} \bar{X}_{1}-\sqrt{\frac{B}{AD}}w_{4}
\bar{X}_{2},
\quad [W,\bar{X}_{2}]= \frac{1}{2 \sqrt{D}} w_{4} \bar{X}_{2} ,
\quad [W,\bar{X}_{3}]= - \frac{1}{\sqrt{D}} w_{4} \bar{X}_{3},\\
[W,\bar{X}_{4}]
& = -\frac{1}{2\sqrt{D}} w_{1} \bar{X}_{1} +\left(\sqrt{\frac{B}{AD}} w_{1} 
- \frac{1}{2 \sqrt{D}} w_{2} \right) \bar{X}_{2} 
+ \frac{1}{\sqrt{D}} w_{3} \bar{X}_{3}.
\end{align*}
We have from (\ref{eq 1})
\begin{align*}
\operatorname{Ric}(W,W)= - \frac{B}{2 AD}w_{1}^{2} + \frac{B}{2AD}
w_{2}^{2} + 0 \cdot w_{3}^{2}- \frac{1}{2}( \frac{3}{D} +\frac{B}{AD}
)w_{4}^{2},
\end{align*}
so the Ricci flow is
\begin{align*}
&  \frac{dA}{dt}= \frac{B}{D}=\frac{B}{AD}A, & & \frac{dB}{dt}= -
\frac{B^{2}}{AD}=-\frac{AD}{B}B,\\
&  \frac{dC}{dt}=0, & &  \frac{dD}{dt}= 3+ \frac{B}{A}.
\end{align*}

It is clear that
\[
C(t)= \lambda_3. \tag{4} \label{eq 4} 
\]
From  $ \frac{1}{A}\frac{dA}{dt}+\frac{1}{B}\frac{dB}{dt}=0 $, we get
$AB=\lambda_1 \lambda_2$. Since $ \frac{d}{dt}\left[
\frac{A}{B}\right]  =\frac{2}{D}$,
$A/B$ is increasing, so
\[
3\leq\frac{dD}{dt}=3+\frac{B}{A}\leq3+\frac{\lambda_{2}}{\lambda_{1}},
\]
from which we get
\[
3t+\lambda_{4}\leq D(t)
\leq(3+\frac{\lambda_{2}}{\lambda_{1}})t+\lambda_{4}.
\tag{5} \label{eq 5}
\]

From
\[
\frac{d}{dt}\left[  \frac{AD}{B}\right]
=3(1+\frac{A}{B})\geq3(1+\frac{\lambda_1
}{\lambda_2})\doteq k_{1} \geq 3,
\]
 we get by integrating
\[
\frac{AD}{B}\geq k_{1}t+k_{2}
\]
where $k_2\doteq \frac{\lambda_1 \lambda_4}{\lambda_2}$. 
Using this, we have
\[
\frac{1}{A}\frac{dA}{dt}=\frac{B}{AD}\leq\frac{1}{k_{1}t+k_{2}}
\]
which gives
\[
A(t)\leq k_{3}(k_{1}t+k_{2})^{1/k_{1}}
\]
where $k_3 \doteq \lambda_3 k_2^{-1/k_1}$. 
To get a lower bound of $A(t)$, we compute
\[
\frac{1}{\lambda_1 \lambda_2}\frac{dA^2}{dt}
=\frac{d}{dt}(\frac{A^2}{AB}) 
=\frac{d}{dt}(\frac{A}{B}) =
\frac{2}{D} \geq
\frac{2\lambda_1}{(3\lambda_1+\lambda_2)t+\lambda_1\lambda_4}
\]
 and get by integrating
 \[
 \lambda_1 \sqrt{\frac{2\lambda_2 }{3\lambda_1+\lambda_2} \ln \left(
 1+\frac{3\lambda_1+\lambda_2}{\lambda_1\lambda_4}t \right) +1 } \leq
A(t) \leq k_{3}(k_{1}t+k_{2})^{1/k_{1}} \tag{6} \label{eq 6}
 \]

From $AB=\lambda_1 \lambda_2$ we get
\[
 \frac{\lambda_1 \lambda_2}{k_{3}}(k_{1}t+k_{2})^{-1/k_1}\leq B(t)
 \leq \frac{\lambda_2 }{\sqrt{\frac{2\lambda_2
}{3\lambda_1+\lambda_2} \ln \left(
 1+\frac{3\lambda_1+\lambda_2}{\lambda_1\lambda_4}t \right) +1 }}. 
 \tag{7} \label{eq 7}
\]
Hence for the family in Proposition 4, the long time behavior of the
solution $g(t)$ as $t \to +\infty$ is
\[
A(t) \to +\infty, \qquad B(t) \to 0^+, \qquad C(t)= \lambda_3, 
\qquad D(t) \to + \infty.
\]
If $k_1>4$, then $(M_q,g_N(t),p)$ converges/collapses to a line in
the pointed Gromov-Hausdorff topology.

Next we compute the curvature decay of $g(t)$. From (\ref{eq 2}) we
find
\begin{xalignat*}{3}
  U(X_1,X_1) & =\frac{A}{2D}X_4
& U(X_2,X_2) & =\frac{B}{2D}X_4
& U(X_3,X_3) & =-\frac{C}{D}X_4   \\
  U(X_4,X_4) & =0
& U(X_1,X_2) & =-\frac{B}{2D}X_4
& U(X_1,X_3) & =0 \\
  U(X_2,X_3) & =0
& U(X_1,X_4) & = -\frac{1}{4}X_1
& U(X_2,X_4) & =\frac{B}{2A}X_1-\frac{1}{4}X_2 \\
  U(X_3,X_4) & =\frac{1}{2}X_3.  
\end{xalignat*} 
From (\ref{eq 3}) with $h=g$ we find the sectional curvatures
\begin{xalignat*}{3}
  K(X_1,X_2) & =\frac{-1 + \frac{B}{A}}{4D} 
& K(X_1,X_3) & =\frac{1}{2D}
& K(X_2,X_3) & =\frac{1}{2D}   \\
  K(X_1,X_4) & = -\frac{1+3\frac{B}{A}}{4D} 
& K(X_2,X_4) & = \frac{-1+\frac{B}{A}}{4D}
& K(X_3,X_4) & =-\frac{1}{D}.
\end{xalignat*} 
Hence for the family in Proposition 4, the curvatures of the solution
$g(t)$ decay at the rate $1/t$.

\subsection*{A6. $U1[3]$}

For $U1[3]$ we use $Y_{i}=\Lambda^{k}_{\phantom{k} i}X_{k}$ with
\[
\Lambda=
\begin{bmatrix}
1 & a_{2} & a_{3} & 0\\
0 & 1 & a_{1} & 0\\
0 & 0 & 1 & 0\\
a_{4} & a_{5} & a_{6} & 1
\end{bmatrix}
\]
to diagonalize the initial metric $g_0$.

\begin{proposition}
For the class $U1[3]$ suppose the initial metric $g_{0}$ is diagonal
in the basis $Y_i$.  Then the Ricci flow solution $g(t)$ remains
diagonal in the basis $Y_i$ if and only if $a_1=a_2$.
\end{proposition}

\begin{proof} We compute
\begin{align*}
[Y_{2},Y_{3}]  &  =0, & [Y_{3},Y_{1}]  &  =0, & [Y_{1},Y_{2}]  &
=0,\\
[Y_{1},Y_{4}]  &  =Y_{2}+\alpha Y_{2}, & [Y_{2},Y_{4}]  &  =Y_{3}, &
[Y_{3},Y_{4}]  &  =0.
\end{align*}
where $\alpha=a_{2}-a_{1}$. We compute the off-diagonal components of 
the Ricci tensor in the basis $\{\bar{Y}_{i}\}$ as in \S 2.A2 and get
\[
\operatorname{Ric}(\bar{Y}_{1},\bar{Y}_{2})= -\frac{\alpha
\lambda_{3}}{2\sqrt{\lambda_{1}\lambda_{2}}\lambda_{4}}, \qquad
 \operatorname{Ric}(\bar{Y}_{2},\bar{Y}_{3})= \frac{\alpha
\sqrt{\lambda_{2}\lambda_{3}}} {2\lambda_{1}\lambda_{4}},
\]
and $ \operatorname{Ric}(\bar{Y}_{1},\bar{Y}_{3})=
\operatorname{Ric}(\bar{Y}_{1},\bar{Y}_{4})
=\operatorname{Ric}(\bar{Y}_{2},\bar{Y}_{4}) =
\operatorname{Ric}(\bar{Y}_{3},\bar{Y}_{4}) =0$.  In order for these
off-diagonal components to be zero, we must have $\alpha=0$ and the
proposition is proved.
\end{proof}

If $\alpha=0$, the $Y_i$ and $X_i$ both satisfy the same Lie bracket
relations.  As in \S 2.A2, we use $X_i$ and carry out the analysis of
the long time behavior of the Ricci flow solutions for the family in
Proposition 5.  Proceeding as in \S 2.A2, we compute
\begin{align*}
&  [W,\bar{X}_{1}]= - \sqrt{\frac{B}{AD}} w_{4} \bar{X}_{2} &&
[W,\bar{X}_{2}]= - \sqrt{\frac{C}{BD}} w_{4} \bar{X}_{3},\\
&  [W,\bar{X}_{3}]= 0 &&
[W,\bar{X}_{4}]= \sqrt{\frac{B}{AD}} w_{1} \bar{X}_{2} + \sqrt
{\frac{C}{BD}} w_{2} \bar{X}_{3}.
\end{align*}
We have from (\ref{eq 1})
\[
\operatorname{Ric}(W,W)= - \frac{B}{2 AD}w_{1}^{2} +
\frac{1}{2}(\frac{B}{AD} - \frac{C}{B D}) w_{2}^{2} + \frac{C}{2 BD}
w_{3}^{2} 
- \frac{1}{2}( \frac{B}{AD} +\frac{C}{BD} )w_{4}^{2},
\]
so the Ricci flow is
\begin{align*}
&  \frac{dA}{dt}= \frac{B}{D} &&  \frac{dB}{dt}= \frac{AC- B^{2}}{AD}
\\
&  \frac{dC}{dt}= - \frac{C^{2}}{ BD}
& &  \frac{dD}{dt}= \frac{B}{A}+ \frac{C}{B}.
\end{align*}

Note that
\begin{align*}
&  \frac{1}{A}\frac{dA}{dt}= \frac{B}{AD} && 
\frac{1}{B}\frac{dB}{dt}= \frac{C}{BD}- \frac{B}{AD} \\
&  \frac{1}{C}\frac{dC}{dt}= - \frac{C}{BD} 
&&  \frac{1}{D}\frac{dD}{dt}= \frac{B}{AD}+ \frac{C}{BD}.
\end{align*}
Hence $ABC=\lambda_1 \lambda_2 \lambda_3$ 
and $\frac{A}{CD}=\frac{\lambda_1}{\lambda_3 \lambda_4}$.
Define $E \doteq \frac{B}{AD}$ and $F \doteq \frac{C}{BD}$, we
compute
\[ 
\frac{dE}{dt}= -3E^{2} \qquad  \frac{dF}{dt}= -3F^{2}.
\]
Solving these gives
\[
E(t)=\frac{E_{0}}{3E_{0}t+1} \qquad F(t)=\frac{F_{0}}{3F_{0}
t+1}
\]
where $E_{0}\doteq \frac{\lambda_2}{\lambda_1\lambda_4}$ and
$F_{0}\doteq \frac{\lambda_3}{\lambda_2 \lambda_4}$.  Using these in
the equations for $(1/A)(dA/dt)$ and $(1/C)(dC/dt)$, we can integrate
to get
\[
A(t)=\lambda_1 \left( 3E_{0}t+1 \right)^{1/3} \qquad 
C(t)=\lambda_3 \left(  3F_{0}t+1 \right)
^{-1/3}. \tag{8} \label{eq 8}
\]
Using these with the conserved quantities $ABC$ and $\frac{A}{CD}$,
we get
\begin{align*}
B(t)  &  = \lambda_{2}\left(  3E_{0}t+1\right)^{-1/3} \left(
3F_{0}t+1\right)  ^{1/3} \\
D(t)  &  =\lambda_{4}\left(  3E_{0}t+1\right)^{1/3} \left(
3F_{0}t+1\right)^{1/3}. \tag{9} \label{eq 9}
\end{align*}
Hence for the family in Proposition 5, the long time behavior of the
solution $g(t)$ is
\[ 
A(t) \to +\infty \qquad B(t) \to (
\lambda_{1}\lambda_{2}\lambda_{3})^{1/3}
\qquad C(t) \to 0^+ \qquad D(t) \to +\infty.
\]
The volume-normalized solution $(M_q,g_N(t),p)$ converges/collapses
to a plane in the pointed Gromov-Hausdorff topology.

Next we compute the curvature decay of $g(t)$. From (\ref{eq 2}) we
find
\begin{xalignat*}{3}
  U(X_1,X_1) & =0
& U(X_2,X_2) & = 0
& U(X_3,X_3) & = 0   
\\
  U(X_4,X_4) & =0
& U(X_1,X_2) & = -\frac{B}{2D}X_4
& U(X_1,X_3) & =0
\\
  U(X_2,X_3) & =-\frac{C}{2D}X_4
& U(X_1,X_4) & = 0 
& U(X_2,X_4) & = \frac{B}{2A}X_1 \\
  U(X_3,X_4)& =\frac{C}{2B}X_2. 
\end{xalignat*} 
From (\ref{eq 3}) with $h=g$ we find the sectional
curvatures
\begin{xalignat*}{3}
  K(X_1,X_2) & =\frac{ \frac{B}{A}}{4D} 
& K(X_1,X_3) & =0
& K(X_2,X_3) & =\frac{\frac{C}{B}}{4D}   \\
  K(X_1,X_4) & = -\frac{3\frac{B}{A}}{4D} 
& K(X_2,X_4) & = \frac{\frac{B}{A}-3\frac{C}{B}}{4D}
& K(X_3,X_4) &=\frac{\frac{C}{B}}{4D}.
\end{xalignat*} 
Hence for the family in Proposition 5, the curvatures of the solution
$g(t)$ decay at the rate $1/t$.

\subsection*{A7. $U3I0$}

For $U1[2,1]$, we use $Y_{i}=\Lambda^{k}_{\phantom{k} i}X_{k}$ with
\[
\Lambda=
\begin{bmatrix}
1 & a_{4} & a_{5} & a_{6}\\
0 & 1 & a_{2} & a_{3}\\
0 & 0 & 1 & a_{1}\\
0 & 0 & 0 & 1
\end{bmatrix}
\]
to diagonalize the initial metric $g_0$.

\begin{proposition}
For the class $U3I0$ suppose the initial metric $g_{0}$ is diagonal in
the basis $Y_i$.  Let
\[
\alpha \doteq a_2 \qquad \beta \doteq a_{1}a_{2}-a_{3}-a_{4}
\qquad
\gamma \doteq a_{1}-a_{1}a_{2}^{2}+a_{2}a_{3}+a_{2}a_{4}-a_{5}.
\]
Then the Ricci flow solution $g(t)$ remains diagonal in the basis
$Y_{i}$ if and only if either

(i) $\alpha=\beta=\gamma=0$; or

(ii) $\beta=\gamma=0$ and $\lambda_2 = (1- \alpha^2 )\lambda_3$.
\end{proposition}

\begin{proof} We compute
\begin{align*}
&  [Y_{1},Y_{4}]=0 \qquad[Y_{2},Y_{4}]=0 && [Y_{3},Y_{4}]=0 \qquad
[Y_{2},Y_{3}]=Y_{4}\\
&  [Y_{3},Y_{1}]=Y_{2}-\alpha Y_{3}+\beta Y_{4} &&
[Y_{1},Y_{2}]=-\alpha
Y_{2}+(\alpha^{2}-1)Y_{3}+\gamma Y_{4}.
\end{align*}

We compute the off-diagonal components of the Ricci tensor in the
basis $\{\bar{Y}_{i}\}$ using (\ref{eq 1}) as in \S 2.A2 and get
\begin{xalignat*}{2}
  \operatorname{Ric}(\bar{Y}_{1},\bar{Y}_{2}) 
& = \frac{\beta\lambda_{4}}{2\sqrt{\lambda_{1}\lambda_{2}}\lambda_{3}} 
& \operatorname{Ric}(\bar{Y}_{1},\bar{Y}_{3}) 
& =
\frac{\gamma\lambda_{4}}{2\sqrt{\lambda_{1}\lambda_{3}}\lambda_{2}}\\
 \operatorname{Ric}(\bar{Y}_{1},\bar{Y}_{4}) 
& =0 
&\operatorname{Ric}(\bar{Y}_{2},\bar{Y}_{3}) 
& = \frac{-2\alpha\lambda_{2}+2\alpha(1-\alpha^{2})\lambda_{3}+\beta\gamma\lambda_{4}}
{2\lambda_{1}\sqrt{\lambda_{2}\lambda_{3}}} 
 \\
 \operatorname{Ric}(\bar{Y}_{2},\bar{Y}_{4}) 
& =\frac{(\beta\lambda_{2}-\alpha\gamma\lambda_{3})\sqrt{\lambda_{4}}} {2\lambda_{1}
\sqrt{\lambda_{2}}\lambda_{3}}
& \operatorname{Ric}(\bar{Y}_{3},\bar{Y}_{4}) 
& = \frac{(-\alpha
\beta\lambda_{2}+(\alpha^{2}-1)\gamma\lambda_{3})\sqrt{\lambda_{4}}}
{2\lambda_{1}\lambda_{2}\sqrt{\lambda_{3}}}
\end{xalignat*}
In order for these off-diagonal components to be zero, we must have
either (i) or (ii) in the proposition.  To finish the proof of (ii) we
need to ensure that the condition $B(t) = (1- \alpha^2 ) C(t)$ holds
for all $t>0$.  We prove this near the end of this subsection.
\end{proof}

\vskip .3cm \textbf{A7i}.  First we study the family (i) in
Proposition 6.  If $\alpha=\beta=\gamma=0$, the bases $Y_i$ and $X_i$
both satisfy the same Lie bracket relations.  As in \S 2.A2, we use
$X_i$ and carry out the analysis of the long time behavior of the
Ricci flow solutions for the family (i) in Proposition 6.  Proceeding
as in \S 2.A2, we find
\begin{xalignat*}{2}
  [W,\bar{X}_{1}] & = \sqrt{\frac{B}{AC}} w_{3} \bar{X}_{2} + \sqrt
{\frac{C}{AB}} w_{2} \bar{X}_{3} 
&  [W,\bar{X}_{2}] & = - \sqrt{\frac{C}{AB}} w_{1} \bar{X}_{3} -
\sqrt{\frac{D}{BC}} w_{3} \bar{X}_{4} \\
  [W,\bar{X}_{3}] & = - \sqrt{\frac{B}{AC}} w_{1} \bar{X}_{2} +
\sqrt{\frac{D}{BC}} w_{2} \bar{X}_{4}
&  [W,\bar{X}_{4}] & = 0.
\end{xalignat*}
We have from (\ref{eq 1})
\begin{multline*}
\operatorname{Ric}(W,W) =   - \frac{1}{2}(\frac{B}{AC} +
\frac{C}{AB} + \frac{2}{A})w_{1}^{2}
-\frac{1}{2} (\frac{C}{AB} + \frac{D}{BC} -\frac{B}{AC} )w_{2}^{2}\\
  -\frac{1}{2}( \frac{B}{AC} + \frac{D}{BC} - \frac{C}{AB})
w_{3}^{2} +
\frac{1}{2} \frac{D}{BC} w_{4}^{2},
\end{multline*}
so the Ricci flow equation is
\begin{align*}
&  \frac{dA}{dt}= \frac{B}{C} + \frac{C}{B} + 2
&&  \frac{dB}{dt}= \frac{C}{A} + \frac{D}{C} -\frac{B^{2}}{AC}\\
&  \frac{dC}{dt}=\frac{B}{A} + \frac{D}{B} - \frac{C^{2}}{AB}
&&  \frac{dD}{dt}= -\frac{D^{2}}{ BC}.
\end{align*}

Straightforward calculations give us
\[
\frac{1}{BC}\frac{d}{dt}[BC]=2\frac{D}{BC} \qquad
\frac{1}{D}\frac{dD}{dt}=-\frac{D}{BC},
\]
from which we conclude that $BCD^{2}=\lambda_2 \lambda_3 \lambda_4^2$.
Now we can write $\frac{dD}{dt}=-\frac{D^{4}}{\lambda_2 \lambda_3
\lambda_4^2}$ and solve to get
\[
D(t)=\lambda_4 \left(
1+3\frac{\lambda_4}{\lambda_2\lambda_3}t\right)^{-1/3}.
\tag{10} \label{eq 10}
\] 

Another set of simple calculations give us
\begin{align*}
  \frac{1}{B-C}\frac{d}{dt}[B-C]=\frac{AD-(B+C)^{2}}{ABC} \qquad
 \frac{1}{AD}\frac{d}{dt}[AD]=-\frac{AD-(B+C)^{2}}{ABC}
\end{align*}
from which we conclude that $AD(B-C)=\lambda_1 \lambda_4(\lambda_2
-\lambda_3)$.  We get
\[
\frac{(B-C)^{2}}{BC}A^{2} =
\frac{A^{2}D^{2}(B-C)^{2}}{BCD^{2}}
=\frac{\lambda_1^2 (\lambda_2-\lambda_3)^2}{\lambda_2
\lambda_3}\doteq 4k_3^2.
\]
for some $k_3 \geq 0$. With this and the identity
\[
\frac{(B+C)^{2}}{BC}=\frac{(B-C)^{2}}{BC}+4
\]
we can write
\[
\frac{dA}{dt}= \frac{4 k_3^2}{A^{2}}+4.
\]
Integrating gives us
\[
A - {k_3} \tan^{-1}\left(  \frac
{A}{k_3}\right)  = 4t+\lambda_{1}- {k_3} \tan^{-1}\left(  \frac
{\lambda_1}{k_3}\right). \tag{11} \label{eq 11}
\]
For large $t$, we have $A(t)\sim 4t$.

Using the conserved quantities $AD(B-C)$ and $BCD^{2}$, we could solve
to get $B(t)$ and $C(t)$ explicitly.  More importantly, we see that
for large $t$,
\[
B(t)\sim C(t)\sim\frac{1}{D(t)}\sim t^{1/3}. \tag{12} \label{eq 12}
\]
Hence for the family (i) in Proposition 6, the long time behavior of
the solution $g(t)$ is
\[ 
A(t) \to +\infty \qquad B(t) \to +\infty  
\qquad C(t) \to +\infty  \qquad D(t) \to 0^+.
\]
For the volume-normalized Ricci flow, the metric components have the
following long time behavior: $A_N(t)\to+\infty$, $B_N(t)$ and
$C_N(t)$ approach some positive constant, and $D_N(t)\to 0^+$.

Next we compute the curvature decay of $g(t)$. From (\ref{eq 2}) we
find
\begin{xalignat*}{3}
  U(X_1,X_1) & =0
& U(X_2,X_2) & =0
& U(X_3,X_3) & = 0   \\
  U(X_4,X_4) & =0
& U(X_1,X_2) & =\frac{B}{2C}X_3  
& U(X_1,X_3) & =\frac{C}{2B}X_2 \\
  U(X_2,X_3) & =-\frac{B+C}{2A}X_1 
& U(X_1,X_4) & = 0
& U(X_2,X_4) & = -\frac{D}{2C}X_3 \\
  U(X_3,X_4) & =\frac{D}{2B}X_2.  
\end{xalignat*} 
From (\ref{eq 3}) with $h=g$ we find the sectional curvatures
\begin{xalignat*}{3}
  K(X_1,X_2) & =\frac{\frac{B}{C}-3 \frac{C}{B}-2}{4A} 
& K(X_1,X_3) & =\frac{\frac{C}{B}-3\frac{B}{C}-2}{4A}
& K(X_1,X_4) & = 0 
\\
  K(X_2,X_3) & =-\frac{3D}{4BC} +\frac{\frac{B}{C}+\frac{C}{B}+2}{4A} 
& K(X_2,X_4) & = \frac{D}{4BC}
& K(X_3,X_4) & =\frac{D}{4BC}.
\end{xalignat*} 
Hence for the family (i) in Proposition 6, the curvatures of the
solution $g(t)$ decay at the rate $1/t$.

\vskip .3cm \textbf{A7ii}.  For the rest of this subsection, we
address the family (ii) in Proposition 6.  Suppose $\beta=\gamma=0$.
The Lie brackets from the proof of Proposition 6 take the form
\begin{align*}
&  [Y_{1},Y_{4}]=0 \qquad[Y_{2},Y_{4}]=0 && [Y_{3},Y_{4}]=0 \qquad
[Y_{2},Y_{3}]=Y_{4}\\
&  [Y_{3},Y_{1}]=Y_{2}-\alpha Y_{3} && [Y_{1},Y_{2}]=-\alpha
Y_{2}+(\alpha^{2}-1)Y_{3}.
\end{align*}
Recall we must show the condition $B(t)=(1-\alpha^2)C(t)$ is preserved
under Ricci flow.  Let $\omega_i$ be the frame dual to $Y_i$.  Assume
the Ricci flow solution $g$ takes the form
\[
g(t)=
A(t)(\omega_1)^2+B(t)(\omega_2)^2+C(t)(\omega_3)^2+D(t)(\omega_4)^2
\]
with
\[
g_0= \lambda_1(\omega_1)^2+
\lambda_2(\omega_2)^2+\lambda_3(\omega_3)^2
+ \lambda_4 (\omega_4)^2.
\]
 
Let $\bar{Y}_{1} \doteq \frac{1}{\sqrt{A}}Y_1,\cdots, \bar{Y}_{4}
\doteq \frac{1}{\sqrt{D}}Y_4$ and let
$W=w_1\bar{Y}_{1}+w_2\bar{Y}_{2}+w_3\bar{Y}_{3}+w_4\bar{Y}_{4}$.  We
first compute $[W,\bar{Y}_i]$ and then compute the Ricci curvature of
$g(t)$ using (\ref{eq 1}).  We find that $\operatorname{Ric}(W,W)$ is
given by
\begin{multline*}
    \operatorname{Ric}(W,W)= -\frac{B^2+2(1+\alpha^2)BC+(1-\alpha^2)^2C^2}{2ABC}w_1^2
+ \frac{-AD+B^2-(1-\alpha^2)^2C^2}{2ABC}w_2^2\\
 +\frac{-AD-B^2+(1-\alpha^2)^2C^2}{2ABC}w_3^2 
 +\frac{D}{2BC}w_4^2 + \frac{\alpha \left(-B+(1-\alpha^2)C\right)
 }{A\sqrt{BC}}w_2w_3.
\end{multline*}
The Ricci flow equation is 
\begin{align*}
& \frac{dA}{dt}= \frac{B^2+2(1+\alpha^2)BC+(1-\alpha^2)^2C^2}{BC}
&& \frac{dB}{dt}= \frac{AD -B^2+ (1-\alpha^2)^2C^2}{AC} \\
& \frac{dC}{dt}= \frac{AD+ B^2-(1-\alpha^2)^2C^2}{AB}
&& \frac{dD}{dt}=-\frac{D^2}{BC}. 
\end{align*}
Hence
\begin{align*}
\frac{d}{dt}(-B+(1-\alpha^2)C)=\frac{AD-(B+(1-\alpha^2)C)^2}{ABC}
\left(-B+(1-\alpha^2)C \right).
\end{align*}
Since $-B+(1-\alpha^2)C=0$ at time $t=0$, it remains $0$ for all time
which implies that $g(t)$ is diagonal in the basis $Y_i$.

With $-B+(1-\alpha^2)C=0$, the Ricci flow equations reduce to
\begin{align*}
& \frac{dA}{dt}= 4
&& \frac{dB}{dt}= (1-\alpha^2)\frac{D}{B} 
&& \frac{dD}{dt}=-(1-\alpha^2)\frac{D^2}{B^2}. 
\end{align*}
A simple calculation gives $\frac{d}{dt}(B^2D^2)=0$ which implies
$BD=\lambda_2 \lambda_4$.  Now we can solve the Ricci flow equations
to obtain
\begin{align*}
& A=\lambda_1+ 4t 
&& B=\left( \lambda_2^3 + 3(1-\alpha^2)\lambda_2 \lambda_4 t
\right)^{1/3} \\
& C=\frac{1}{1-\alpha^2} \left( \lambda_2^3 + 3(1-\alpha^2)\lambda_2 
\lambda_4 t \right)^{1/3} 
&&D={\lambda_2 \lambda_4} \left( \lambda_2^3 + 3(1-\alpha^2)\lambda_2 
\lambda_4 t \right)^{-1/3}. 
\end{align*}

Hence for the family (ii) in Proposition 6, the long time behavior of
the Ricci flow $g(t)$ as $t \to \infty$ is
\[
A(t) \to +\infty \qquad B(t) \to  +\infty \qquad C(t) \to  +\infty 
\qquad D(t) \to  0^+
\]  
For the volume-normalized Ricci flow, the metric components have the
following long time behavior: $A_N(t)\to+\infty$, $B_N(t)$ and
$C_N(t)$ approach some positive constants, and $D_N(t)\to 0^+$.

Finally we compute the curvature decay of $g(t)$.  From (\ref{eq 2})
we find
\begin{xalignat*}{3}
  U(Y_1,Y_1) & =0 
& U(Y_2,Y_2) & =-\frac{\alpha B}{A}Y_1
& U(Y_3,Y_3) & = \frac{\alpha C}{A}Y_1   \\
  U(Y_4,Y_4) & =0 
& U(Y_1,Y_2) & =\frac{\alpha}{2}Y_2+\frac{B}{2C}Y_3
& U(Y_1,Y_3) & = \frac{(1-\alpha^2)C}{2B}Y_2-\frac{\alpha}{2}Y_3  \\
U(Y_1,Y_4) & = 0 
& U(Y_2,Y_3) & =-\frac{B+(1-\alpha^2)C}{2A}Y_1 
& U(Y_2,Y_4) & = -\frac{D}{2C}Y_3 \\
  U(Y_3,Y_4) & =\frac{D}{2B}Y_2.  
\end{xalignat*} 
From (\ref{eq 3}) with $h=g$ we find the sectional curvatures
\begin{xalignat*}{3}
  K(Y_1,Y_2) & =-\frac{1}{A} 
& K(Y_1,Y_3) & =-\frac{1}{A} 
& K(Y_2,Y_3) & =-\frac{3D}{4BC}+\frac{1}{A}  \\
  K(Y_1,Y_4) & = 0
& K(Y_2,Y_4) & = \frac{D}{4BC} 
& K(Y_3,Y_4) & =\frac{D}{4BC}.
\end{xalignat*} 
Hence for the family (ii) in Proposition 6, the curvatures of the
solution $g(t)$ decay at the rate $1/t$.

\subsection*{A8. $U3I2$}

For $U3I2$, we use $Y_{i}=\Lambda^{k}_{\phantom{k} i}X_{k}$ with
\[
\Lambda=
\begin{bmatrix}
1 & a_{4} & a_{5} & a_{6}\\
0 & 1 & a_{2} & a_{3}\\
0 & 0 & 1 & a_{1}\\
0 & 0 & 0 & 1
\end{bmatrix}
\]
to diagonalize the initial metric $g_0$.

\begin{proposition}
For the class $U3I2$ suppose the initial metric $g_{0}$ is diagonal in
the basis $Y_i$.  Then the Ricci flow solution $g(t)$ remains diagonal
if and only if $a_2=0$, $a_1=a_5$ and $a_3=a_4$.
\end{proposition}

\begin{proof} We compute
\begin{align*}
&  [Y_{1},Y_{4}]=0 \qquad[Y_{2},Y_{4}]=0 && [Y_{3},Y_{4}]=0 \qquad
[Y_{2},Y_{3}]=-Y_{4}\\
&  [Y_{3},Y_{1}]=Y_{2}+\alpha Y_{3}+\beta Y_{4} &&
[Y_{1},Y_{2}]=\alpha
Y_{2}+(1+\alpha^{2})Y_{3}+\gamma Y_{4}.
\end{align*}
where
\[
\alpha=-a_{2} \qquad \beta=a_{1}a_{2}-a_{3}+a_{4} \qquad
\gamma=-a_{1}-a_{1}a_{2}^{2}+a_{2}a_{3}-a_{2}a_{4}+a_{5}.
\]
We compute the off-diagonal components of the Ricci tensor 
in the basis $\{\bar{Y}_{i}\}$ using (\ref{eq 1}) as in \S 2.A2 and
get
\begin{xalignat*}{2}
\operatorname{Ric}(\bar{Y}_{1},\bar{Y}_{2}) 
& = -\frac{\beta\lambda_{4}}{2\sqrt{\lambda_{1}\lambda_{2}}\lambda_{3}} 
& \operatorname{Ric}(\bar{Y}_{1},\bar{Y}_{3})
& = -\frac{\gamma\lambda_{4}}{2\sqrt{\lambda_{1}\lambda_{3}}\lambda_{2}}\\
  \operatorname{Ric}(\bar{Y}_{1},\bar{Y}_{4}) 
& =0 
& \operatorname{Ric}(\bar{Y}_{2},\bar{Y}_{3}) 
& = \frac{2\alpha \lambda_{2}+2\alpha(1+\alpha^{2})\lambda_{3}+\beta\gamma\lambda_{4}}
{2\lambda_{1}\sqrt{\lambda_{2}\lambda_{3}}} \\
\operatorname{Ric}(\bar{Y}_{2},\bar{Y}_{4})
& = \frac{(\beta\lambda_{2}+\alpha\gamma\lambda_{3})\sqrt{\lambda_{4}}} {2\lambda_{1}
\sqrt{\lambda_{2}}\lambda_{3}}
& \operatorname{Ric}(\bar{Y}_{3},\bar{Y}_{4})
& = \frac{(\alpha\beta\lambda_{2}+(1+\alpha^{2})\gamma\lambda_{3})\sqrt{\lambda_{4}}}
{2\lambda_{1}\lambda_{2}\sqrt{\lambda_{3}}}.
\end{xalignat*}
In order for these off-diagonal components to be zero, we must have
$\alpha=\beta=\gamma=0$ and the proposition follows.
\end{proof}

If $\alpha=\beta=\gamma=0$, the bases $Y_i$ and $X_i$ both satisfy the
same Lie bracket relations.  As in \S 2.A2, we use $X_i$ and carry out
the analysis of the long time behavior of the Ricci flow solutions for
the family in Proposition 7.  Proceeding as in \S 2.A2, we find
\begin{align*}
&  [W,\bar{X}_{1}]= \sqrt{\frac{B}{AC}} w_{3} \bar{X}_{2} -
\sqrt{\frac{C}
{AB}} w_{2} \bar{X}_{3}
& &  [W,\bar{X}_{2}]= \sqrt{\frac{C}{AB}} w_{1} \bar{X}_{3} +
\sqrt{\frac{D}
{BC}} w_{3} \bar{X}_{4}\\
&  [W,\bar{X}_{3}]= - \sqrt{\frac{B}{AC}} w_{1} \bar{X}_{2}
-\sqrt{\frac
{D}{BC}} w_{2} \bar{X}_{4}
& &  [W,\bar{X}_{4}]= 0.
\end{align*}
We have from (\ref{eq 1})
\begin{multline*}
\operatorname{Ric}(W,W)=
\frac{1}{2}(\frac{2}{A}- \frac{C}{AB}-\frac{B}{AC} ) w_{1}^{2} 
+ \frac{1}{2}(\frac{B}{AC}- \frac{C}{AB}-\frac{D}{BC} ) w_{2}^{2}\\
  +\frac{1}{2}(\frac{C}{AB}- \frac{B}{AC}- \frac{D}{BC} )w_{3}^{2} +
\frac{D}{2BC} w_{4}^{2}.
\end{multline*}
The Ricci flow is
\begin{align*}
&  \frac{dA}{dt}= \frac{C}{B}+ \frac{B}{ C} -2 
&&  \frac{dB}{dt}= - \frac{B^{2}}{AC} + \frac{C}{A} + \frac{D}{C} \\
&  \frac{dC}{dt}= - \frac{C^{2}}{AB} + \frac{B}{A}+ \frac{D}{B}
&&  \frac{dD}{dt}= - \frac{D^{2}}{BC}.
\end{align*}

The equations here are similar to those of the case A7(i), with the
only difference being in the equation for $A$.  Because the equations
for $B$, $C$ and $D$ are the same, we know that $BCD^{2}=\lambda_2
\lambda_3 \lambda_4^2$.  It follows that $ \frac{dD}{dt}= -
\frac{D^{4}}{\lambda_2 \lambda_3 \lambda_4^2}$, and hence
\[
D(t)=\lambda_{4}\left(  1+ \frac{3\lambda_{4}}{\lambda_{2}
\lambda_{3}}t\right)  ^{-1/3}. \tag{13} \label{eq 13}
\]

Calculations similar to those in the case A7(i) show that $AD(B+C)=
\lambda_1 \lambda_4 (\lambda_2+ \lambda_3) $.  So
$A^2(\frac{C}{B}+\frac{B}{C}+2)=\frac{(AD(B+C))^2}{BCD^2} \doteq
k_4^2$ is a constant where $k_4 \geq 0$, and
\[
\frac{dA}{dt}=\frac{k_4^2}{A^{2}}-4.
\]
Integrating the equation gives us
\[
\frac{k_4}{2} \tanh^{-1}\left(  \frac{2A}{k_4}\right)  
- A= 4t+ k_5 \tag{14} \label{eq 14}
\]
where $k_5$ is a constant.  Since $A$ increases for all $t$ and $\tanh
x$ asymptotes to $1$, we see that $A(t)\to k_4/2$ as $t\to +\infty$.

Using the conserved quantity $BCD^2$ and $AD(B+C)$, we conclude that
for the family in Proposition 7 both $B$ and $C$ grow at the rate
$t^{1/3}$, and the long time behavior of the Ricci flow $g(t)$ as $t
\to +\infty$ is
\[
A(t) \to k_4/2 \qquad B(t) \to +\infty \qquad C(t) \to +\infty \qquad
D(t) \to 0^+.
\]  
The volume-normalized flow $(M_q,g_N(t),p)$ converges/collapses to a
plane in the pointed Gromov-Hausdorff topology.

Next we compute the curvature decay of $g(t)$. From (\ref{eq 2}) we
find
\begin{align*}
&U(X_1,X_1)= 0  && U(X_2,X_2)= 0 && U(X_3,X_3)= 0   \\
& U(X_4,X_4)=0  && U(X_1,X_2)= \frac{B}{2C}X_3 && 
U(X_1,X_3)= -\frac{C}{2B}X_2 \\
&U(X_2,X_3)= \frac{-B+C}{2A}X_1 &&  U(X_1,X_4)= 0 && 
 U(X_2,X_4)=\frac{D}{2C}X_3 \\
 & U(X_3,X_4)=-\frac{D}{2B}X_2. && &&  
\end{align*} 
From (\ref{eq 3}) with $h=g$ we find the sectional
curvatures
\begin{align*}
&K(X_1,X_2)=\frac{\frac{B}{C} - 3\frac{C}{B}+2 }{4A} 
&&  K(X_1,X_3)=\frac{ -3 \frac{B}{C}+\frac{C}{B}+2}{4A} 
 \\ 
& K(X_2,X_3)=-\frac{3D}{4BC} + \frac{ \frac{B}{C}+\frac{C}{B}-2}{4A} 
&& K(X_1,X_4)= 0 \\
& K(X_2,X_4)= \frac{D}{4BC} &&
 K(X_3,X_4)=\frac{D}{4BC}.
\end{align*} 

Here the decay rate is not obvious for all sectional curvatures.  Note
that $\frac{B}{C}\to 1$.  The decay rate of $\frac{B}{C}-1$ follows
from the equation of $\frac{dA}{dt}$.  It suffices to show that
$\frac{dA}{dt}$ decays at the rate $e^{-ct}$ for some $c>0$.  From
(\ref{eq 14}) we get
\[
A= \frac{k_4}{2} \tanh \left[\frac{8}{k_4}t +
\frac{2}{k_4}A+\frac{2k_5}{k_4} \right]
\]
and the decay rate of $\frac{dA}{dt}$ follows from taking the time
derivative of this equation.  Thus $\frac{B}{C}-1$ decays
exponentially.  Hence for the family in Proposition 7, the curvatures
of the solution $g(t)$ decay at the rate $1/t$.

\subsection*{A9. $U3S1$}

The Lie bracket relations for cases A9 and A10 differ only in
$[X_1,X_2]$.  To unify some of the calculations for these two cases,
we introduce a constant $\delta$ and write $[X_1,X_2] = \delta X_3 $
with $\delta = -1$ corresponding to A9 and $\delta = 1$ corresponding
to A10.

For $U3S1$ and $U3S3$, we use $Y_{i}=\Lambda^{k}_{\phantom{k} i}X_{k}$ with
\[
\Lambda=
\begin{bmatrix}
1 & 0 & 0 & 0 \\
0 & 1 & 0 & 0 \\
0 & 0 & 1 & 0 \\
a_1 & a_2 & a_3 & 1
\end{bmatrix}
\]
to diagonalize the initial metric $g_0$.


\begin{proposition}
For the class $U3S1$ suppose the initial metric $g_{0}$ is diagonal in
the basis $Y_i$.  Then

(i) if $\lambda_1 \neq \lambda_2$, 
the Ricci flow solution $g(t)$ remains diagonal if and only if 
$a_1=a_2=a_3=0$; and

(ii) if  $\lambda_1=\lambda_2$, 
 the Ricci flow solution $g(t)$ remains diagonal if and only if
$a_1=a_2=0$.
\end{proposition}

\begin{proof} We compute
\begin{align*}
[Y_{1},Y_{4}]  &  =-a_{3}Y_{2}+\delta a_{2}Y_{3} & [Y_{2},Y_{4}]  &
= a_{3}Y_{1} -\delta a_{1}Y_{3} & [Y_{3},Y_{4}]  &  = -a_{2}
Y_{1}+ a_{1}Y_{2} \\
[Y_{2},Y_{3}]  &  =Y_{1}, & [Y_{3},Y_{1}]  &  =Y_{2}, & [Y_{1},Y_{2}]
&
=\delta Y_{3}.
\end{align*}
We compute the off-diagonal components of the Ricci tensor in the
basis $ \bar{Y}_{i}$ using (\ref{eq 1}) as in \S 2.A2 and get
\begin{align*}
&  \operatorname{Ric}(\bar{Y}_{1},\bar{Y}_{2})= \frac{(\lambda
_{3}^{2}-\lambda_{1}\lambda_{2}) a_{1} a_{2}} {2\lambda_{3}\lambda
_{4}\sqrt{\lambda_{1}\lambda_{2}}}
&& \operatorname{Ric}(\bar{Y}_{1},\bar{Y}_{3})= \frac{(\lambda
_{2}^{2}-\delta\lambda_{1}\lambda_{3}) a_{1} a_{3}} {2\lambda
_{2}\lambda_{4}\sqrt{\lambda_{1}\lambda_{3}}}\\
&  \operatorname{Ric}(\bar{Y}_{2},\bar{Y}_{3})= \frac{(\lambda
_{1}^{2}-\delta\lambda_{2}\lambda_{3}) a_{2} a_{3}} {2\lambda
_{1}\lambda_{4}\sqrt{\lambda_{2}\lambda_{3}}}
&& \operatorname{Ric}(\bar{Y}_{1},\bar{Y}_{4})= - \frac{(\lambda
_{2}-\delta\lambda_{3})^{2} a_{1}} {2\lambda_{2}\lambda_{3}\sqrt
{\lambda_{1}\lambda_{4}}} \\
&  \operatorname{Ric}(\bar{Y}_{2},\bar{Y}_{4})= - \frac{(\lambda
_{1}-\delta\lambda_{3})^{2} a_{2}} {2\lambda_{1}\lambda_{3}\sqrt
{\lambda_{2}\lambda_{4}}}
&&  \operatorname{Ric}(\bar{Y}_{3},\bar{Y}_{4})= - \frac{(\lambda
_{1}-\lambda_{2})^{2} a_{3}} {2\lambda_{1}\lambda_{2}\sqrt{\lambda
_{3}\lambda_{4}}}.
\end{align*}
The diagonal components are given by 
\begin{align*}
& \operatorname{Ric}(\bar{Y}_{1},\bar{Y}_{1})= \frac{ 
(\lambda_{1}^{2}-\lambda_{3}^{2})\lambda_{2} a_{2}^{2} + (\lambda
_{1}^{2}-\lambda_{2}^{2})\lambda_{3} a_{3}^{2}
+\bigl(\lambda_{1}^{2}-(\lambda
_{2}-\delta\lambda_{3})^{2}\bigr)\lambda_{4}} {2\lambda_{1}\lambda_{2}
\lambda_{3}\lambda_{4}}\\
& \operatorname{Ric}(\bar{Y}_{2},\bar{Y}_{2})= \frac{
(\lambda_{2}^{2}-\lambda_{3}^{2})\lambda_{1}a_{1}^{2} + (\lambda
_{2}^{2}-\lambda_{1}^{2})\lambda_{3}a_{3}^{2}
+\bigl(\lambda_{2}^{2}-(\lambda
_{1}-\delta\lambda_{3})^{2}\bigr)\lambda_{4}} {2\lambda_{1}\lambda_{2}
\lambda_{3}\lambda_{4}}\\
&  \operatorname{Ric}(\bar{Y}_{3},\bar{Y}_{3})= \frac{ 
(\lambda_{3}^{2}-\lambda_{2}^{2})\lambda_{1}a_{1}^{2} + (\lambda
_{3}^{2}-\lambda_{1}^{2})\lambda_{2}a_{2}^{2}
+\bigl(\lambda_{3}^{2}-(\lambda_{1}-\lambda_{2}
)^{2}\bigr)\lambda_{4}} {2\lambda_{1}\lambda_{2}\lambda
_{3}\lambda_{4}}\\
&  \operatorname{Ric}(\bar{Y}_{4},\bar{Y}_{4})= -\frac{ 
(\lambda_{2}-\delta\lambda_{3})^{2}\lambda_{1} a_{1}^{2} +
(\lambda_{1}
-\delta\lambda_{3})^{2} \lambda_{2} a_{2}^{2} +
(\lambda_{1}-\lambda_{2})^{2}\lambda_{3}a_{3}^{2}}
{2\lambda_{1}\lambda_{2}\lambda_{3}\lambda_{4}}. \tag{15} \label{eq
15}
\end{align*}

For $\delta =-1$, in order for these off-diagonal components to be
zero, we must have either (i) or (ii) in Proposition 8.  As in case
A7ii, to finish the proof of (ii) in Proposition 8 we need to ensure
that the condition $A(t)=B(t)$ holds for all $t>0$.  We prove this at
the end of this subsection.
\end{proof}

\noindent \textbf{Remark} Note that there are many initial metrics
$g_{0}$ that cannot be diagonalized by the choice of $\Lambda$ we use
here.  For A9 and A10, the Lie group $G$ is a product
$G_1\times\mathbb{R}$ with $\dim(G_1)=3$.  After transforming with
$\Lambda$ as given above, one can use a Milnor frame on $G_{1}$ (with
respect to a chosen initial metric on $G_{1}$) to further diagonalize,
in which case the Lie algebra takes the form
\begin{align*}
[Y_{1},Y_{4}]  &  =a_{1}Y_{1}+a_{2}Y_{2} +a_{3}Y_{3} \\
[Y_{2},Y_{4}]  &  =b_{1}Y_{1}+b_{2}Y_{2} +b_{3}Y_{3}   \\
[Y_{3},Y_{4}]  &  =c_{1}Y_{1}+c_{2}Y_{2} +c_{3}Y_{3}  \\
[Y_{2},Y_{3}]  &  =Y_{1} 
\qquad [Y_{3},Y_{1}]  =Y_{2} 
\qquad [Y_{1},Y_{2}]   =\delta Y_{3}
\end{align*}
with $a_{1}+b_{2}+c_{3}=0$ from the unimodular condition.
With these, the off-diagonal components of the Ricci curvature are given by
\begin{align*}
\operatorname{Ric}(\bar{Y}_{1},\bar{Y}_{2}) & = \frac
{{c_1}\,{c_2}{{\lambda }_1}{{\lambda
}_2}+b_{1}(b_{2}-a_{1})\lambda_{1}\lambda_{3}+a_{2}(a_{2}-b_{2})\lambda_{2}\lambda_{3}-a_{3}b_{3}\lambda_{3}^{2}}
{2\, {\sqrt{{{\lambda }_1}}}\,{\sqrt{{{\lambda }_2}}}\,
{{\lambda}_3}\,{{\lambda }_4}} \\
\operatorname{Ric}(\bar{Y}_{1},\bar{Y}_{3}) & = 
\frac
{-c_{1}(2a_{1}+b_{2})\lambda_{1}\lambda_{2}+b_{1}b_{3}\lambda_{1}\lambda_{3}
-a_{2}c_{2}\lambda_{2}^{2}+a_{3}(2a_{1}+b_{2})\lambda_{2}\lambda_{3}}
{2\, {\sqrt{{{\lambda }_1}}}\,{{\lambda}_2}\, {\sqrt{{{\lambda }_3}}}\,{{\lambda }_4}} \\
\operatorname{Ric}(\bar{Y}_{1},\bar{Y}_{4}) & = 
-\frac
{ \left( {{\lambda }_2} - \delta \,{{\lambda }_3} \right)\,
\left( {c_2}\,{{\lambda }_2} + {b_3}\,{{\lambda }_3} \right) }
{2\,{\sqrt{{{\lambda }_1}}}\, {{\lambda }_2}\,{{\lambda}_3}\,{\sqrt{{{\lambda}_4}}}}
\\
\operatorname{Ric}(\bar{Y}_{2},\bar{Y}_{3}) & = 
\frac
{-b_{1}c_{1}\lambda_{1}^{2}-c_{2}(a_{1}+2b_{2})\lambda_{1}\lambda_{2}
+b_{3}(a_{1}+2b_{2})\lambda_{1}\lambda_{3}+a_{2}a_{3}\lambda_{2}\lambda_{3}}
{2\,{{\lambda }_1}\, {\sqrt{{{\lambda}_2}}}\,{\sqrt{{{\lambda }_3}}}\, {{\lambda }_4}}
\\
\operatorname{Ric}(\bar{Y}_{2},\bar{Y}_{4}) & = 
\frac
{\left({{\lambda }_1} - \delta \,{{\lambda }_3} \right) \, \left(
{c_1}\,{{\lambda }_1} + {a_3}\,{{\lambda }_3} \right) }
{2\,{{\lambda}_1}\,{\sqrt{{{\lambda }_2}}}\, {{\lambda }_3}\,{\sqrt{{{\lambda}_4}}}}
\\
\operatorname{Ric}(\bar{Y}_{3},\bar{Y}_{4}) & = 
-\frac
{ \left(
{{\lambda }_1} - {{\lambda }_2} \right) \, \left( {b_1}\,{{\lambda
}_1} + {a_2}\,{{\lambda }_2} \right)  }
{2\,{{\lambda}_1}\,{{\lambda }_2}\, {\sqrt{{{\lambda }_3}}}\,{\sqrt{{{\lambda}_4}}}}
\end{align*}
One can analyze these expressions to determine conditions under which 
Ricci flow preserves the diagonalization of an initial metric. 
The complexity of these expressions leads to many cases that must be
analyzed so we have limited our attention to the transformation matrix
$\Lambda$ given above with the results given in Proposition 8 for A9
and Proposition 9 for A10.


\textbf{A9i}.  First we study family (i) in Proposition 8.  For
$a_1=a_2=a_3=0$, we have $Y_i=X_i$.  The metric $g(t)$ is a product
metric on $\widehat{SL}(2,\mathbb{R}) \times \mathbb{R}$
\[
g(t)=g_{SL}(t)+ \lambda^4 du^2
\]
where
$g_{SL}(t)=A(t)(\theta_1)^2+B(t)(\theta_2)^2+C(t)(\theta_3)^2$ 
is a Ricci flow solution on $\widehat{SL}(2,\mathbb{R})$.  
From (\ref{eq 15}), we get the Ricci flow equations
\begin{align*}
& \frac{dA}{dt}=\frac{(B+C)^2-A^2}{BC} 
\quad \frac{dB}{dt}=\frac{(A+C)^2-B^2}{AC} 
\quad \frac{dC}{dt}=\frac{(A-B)^2-C^2}{AB}. 
\end{align*}
The volume-normalized flow associated with $g_{SL}(t)$ has been
analyzed in \cite{IJ}.  It follows that the volume-normalized solution
$(M_q,g_N(t),p)$ converges/collapses to a plane in the pointed
Gromov-Hausdorff topology.  The curvatures of $g(t)$ decay at the rate
$1/t$.

\textbf{A9ii}.  For the rest of this subsection we address family (ii)
in Proposition 8 where $\lambda_1=\lambda_2$ and $a_1=a_2=0$.  From
(\ref{eq 15}) we conclude that the Ricci flow equation of $g(t)$ is
\begin{xalignat*}{2}
  \frac{dA}{dt} & =\frac{(B+C)^2-A^2}{BC}+\frac{-A^2+B^2}{BD}a_3^2 
& \frac{dB}{dt} & =\frac{(A+C)^2-B^2}{AC} +\frac{A^2-B^2}{AD}a_3^2 \\
  \frac{dC}{dt} & =\frac{(A-B)^2-C^2}{AB} 
& \frac{dD}{dt} & = \frac{(A+B)^2}{AB}a_3^2 . 
\end{xalignat*}

Recall we must show that the condition $A(t)=B(t)$ is preserved under
Ricci flow.  To this end, we compute
\[
\frac{d}{dt}(A-B)=\left[\frac{C^2-(A+B)^2}{ABC}-\frac{(A+B)^2
}{ABD}a_3^2\right](A-B).
\]
Since $A-B=0$ at time $t=0$, this implies that $A(t)=B(t)$ and $g(t)$
remains diagonal in the basis $Y_i$.

With $A=B$, the Ricci flow equations reduce to 
\begin{align*}
\frac{dA}{dt}=\frac{C}{A}+2
 \qquad \frac{dC}{dt}=-\frac{C^2}{A^2} 
\qquad \frac{dD}{dt}= 4a_3^2 . 
\end{align*}
A simple computation shows $\frac{d}{dt}\left( \frac{C}{A}\right)
=-2A^{-1} \left( \frac{C}{A} + (\frac{C}{A})^2 \right) \leq 0$; hence
$2 \leq \frac{dA}{dt} \leq 2+\frac{\lambda_3}{\lambda_1}$ and
\begin{align*}
2t+ \lambda_1 \leq A(t)=B(t) \leq (2+\frac{\lambda_3}{\lambda_1})t 
+\lambda_1. \tag{16} \label{eq 16}
\end{align*}
From the equation for $\frac{dC}{dt}$ and (\ref{eq 16}) we get
$-\frac{C^2}{(2t+\lambda_1)^2} \leq \frac{dC}{dt} \leq 0$.
Integrating these inequalities we find
\[
\frac{2\lambda_1 \lambda_3}{2\lambda_1 + \lambda_3} \leq C(t) \leq
\lambda_3. \tag{17} \label{eq 17}
\] 
Finally 
\[D(t)=  4a_3^2 t +\lambda_4. \tag{18} \label{eq 18} 
\]
Hence for family (ii) in Proposition 8 with $a_3\neq 0$, the long time
behavior of of the Ricci flow $g(t)$ as $t \to +\infty$ is
\[
A(t) \to +\infty \qquad B(t) \to +\infty \qquad 
C(t) \to \text{ constant} >0 \qquad D(t) \to +\infty.  
\]

Next we compute the curvature decay of $g(t)$ as in A7ii.  From
(\ref{eq 2}) we find (using $A(t)=B(t)$)
\begin{align*}
& U(Y_1,Y_3)=\frac{A+C}{2A}Y_2 && U(Y_2,Y_3) =-\frac{A+C}{2A}Y_1 \\
& U(Y_1,Y_4) =\frac{1}{2}a_3Y_2 && U(Y_2,Y_4) =-\frac{1}{2}a_3Y_1,
\end{align*} 
and all other $U(Y_i,Y_j)=0$.  From (\ref{eq 3}) with $h=g$ we find
the sectional curvatures
\begin{align*}
K(Y_1,Y_2)=-\frac{4 +3 \frac{C}{A}}{4A} 
\qquad  K(Y_1,Y_3)=K(Y_2,Y_3)=\frac{C}{4A^2}
\end{align*} 
and all other $K(Y_i,Y_j)=0$.  Hence for family (ii) in Proposition
8(ii), the curvatures of the solution $g(t)$ decay at the rate $1/t$.

\subsection*{A10. U3S3}

Using the setup given in A9, we prove the following.
\begin{proposition} 
For the class $U3S1$ suppose the initial metric $g_{0}$ is
diagonal
in the basis $Y_i$. Then 

(i) if $\lambda_1, \lambda_2, \lambda_3$ are all different, the Ricci
flow solution $g(t)$ remains diagonal if and only if $a_1=a_2=a_3=0$;

(ii) if $\lambda_j =\lambda_k \neq \lambda_i$ for some permutation
$\{i,j,k\}$ of $\{1,2,3 \}$, the Ricci flow solution $g(t)$ remains
diagonal if and only if $a_j =a_k=0$; and

(iii) if the initial metric satisfies $\lambda_1=\lambda_2=\lambda_3$,
the Ricci flow solution $g(t)$ remains diagonal for any $a_1$, $a_2$,
and $a_3$.
\end{proposition}

\begin{proof}
Set $\delta =1$ in the proof of Proposition 8.  In order for the
off-diagonal Ricci components to be zero, we must have either (i) or
(ii) or (iii) in Proposition 9.  As in previous cases, to finish the
proof of (ii) in Proposition 9 we need to ensure that the condition
$B(t)=C(t)$ holds for all $t>0$ (using $j=2,k=3$ without loss of
generality).  Also, to finish the proof of (iii) in Proposition 9 we
need to ensure that the condition $A(t)=B(t)=C(t)$ holds for all
$t>0$.  These are verified below.
\end{proof}

\vskip .3cm 
\textbf{A10i}.  
If $a_1=a_2=a_3=0$, we have $Y_i=X_i$.
The metric is a product metric on $S^3 \times \mathbb{R}$
\[
g(t)=g_{S^3}(t)+ \lambda^4 du^2
\]
where $g_{S^3}(t)=A(t)(\theta_1)^2+B(t)(\theta_1)^2+C(t)(\theta_1)^2$
is a Ricci flow solution on $S^3$.  From (\ref{eq 15}), we get the
Ricci flow equations
\[
\frac{dA}{dt}=\frac{(B-C)^2-A^2}{BC} 
\qquad \frac{dB}{dt}=\frac{(A-C)^2-B^2}{AC}
\qquad \frac{dC}{dt}=\frac{(A-B)^2-C^2}{AB}. \tag{19} \label{eq 19}
\]
The volume-normalized flow associated with $g_{S^3}(t)$ has been
analyzed in \cite{IJ} and is found to converge to a round sphere.  It
follows that the volume-normalized solution $(M_q,g_N(t),p)$
converges/collapses to a line in the pointed Gromov-Hausdorff topology.
The curvature behavior of $g(t)$ is a Type I singularity in the sense
of Hamilton.

\vskip .3cm 
\textbf{A10ii}.  
For family (ii) in Proposition 9, without
loss of generality we may assume that $i=1,j=2$, and $k=3$ so
$\lambda_2 =\lambda_3$ and $a_2=a_3=0$.  From (\ref{eq 15}) we
conclude that the Ricci flow equation of $g(t)$ is
\begin{xalignat*}{2}
  \frac{dA}{dt} & =\frac{(B-C)^2-A^2}{BC}  
& \frac{dB}{dt} & =\frac{(A-C)^2-B^2}{AC} -\frac{B^2-C^2}{CD}a_1^2
\\
  \frac{dC}{dt} & =\frac{(A-B)^2-C^2}{AB}+\frac{B^2-C^2}{BD}a_1^2 
& \frac{dD}{dt} & = -\frac{(B-C)^2}{BC}a_1^2 . 
\end{xalignat*}

Recall that we must show that the condition $B(t)=C(t)$ is preserved
under Ricci flow.  This follows from 
\[
\frac{d}{dt}(B-C)=\left[\frac{A^2-(B+C)^2}{ABC}-\frac{(B+C)^2
}{BCD}a_1^2\right](B-C).
\]

With $B=C$, the Ricci flow equations reduce to 
\begin{align*}
\frac{dA}{dt}=-\frac{A^2}{B^2}
 \qquad \frac{dB}{dt}=\frac{A}{B}-2 
\qquad \frac{dD}{dt}= 0 . 
\end{align*}
This is a special case of equation (\ref{eq 19}) with $B=C$, so the
conclusions from A10i hold here.

\vskip .3cm 
\textbf{A10iii}.
For family (iii) in Proposition 9,
$\lambda_1=\lambda_2 =\lambda_3$.  From (\ref{eq 15}) we conclude that
the Ricci flow equation of $g(t)$ is
\begin{align*}
&\frac{dA}{dt}= -\frac{ 
(A^{2}-C^{2})B a_{2}^{2} + (A^{2}-B^{2})C a_{3}^{2} +\bigl(A^{2}-
(B-C)^{2}\bigr)D} {BCD}\\
& \frac{dB}{dt}= -\frac{ (B^{2}-C^{2})Aa_{1}^{2} +
(B^{2}-A^{2})Ca_{3}^{2} +\bigl(B^{2}-(A-C)^{2}\bigr)D} {ACD}\\
&  \frac{dC}{dt}= -\frac{ 
(C^{2}-B^{2})Aa_{1}^{2} + (C^{2}-A^{2})Ba_{2}^{2}
+\bigl(C^{2}-(A-B)^{2}\bigr) D} {ABD}\\
& \frac{dD}{dt} = \frac{ 
(B-C)^{2}A a_{1}^{2} + (A
-C)^{2} B a_{2}^{2} + (A-B)^{2}Ca_{3}^{2}}
{ABC}. 
\end{align*}

Recall we need to show $A(t)=B(t)=C(t)$ is preserved under Ricci
flow.  This follows from
\begin{align*}
&\frac{d}{dt}(A-B)= M_{11}(A-B)+ M_{12}(A-C) \\ 
&\frac{d}{dt}(A-C)= M_{21}(A-B)+ M_{22}(A-C)
\end{align*}
where $M_{ij}$ are continuous functions of $t$.

With $A=B=C$, the Ricci flow equations reduce to 
\begin{align*}
\frac{dA}{dt}= \frac{dB}{dt}=\frac{dC}{dt}=-1  
\qquad \frac{dD}{dt}= 0 , 
\end{align*}
so
\[
g(t)=(\lambda_1-t)(\omega_1)^2+(\lambda_1-t)(\omega_2)^2
+(\lambda_1-t)(\omega_3)^2+\lambda_4 (\omega_4)^2
\] 
where $\omega_i$ is the dual frame of $Y_i$.  It follows from this
explicit solution that the conclusions from A10i hold here.

\section*{3. The Ricci flow of locally homogeneous closed 4-manifolds
modelled
on non-Lie groups}

In this section, all of the metrics are on direct products of spheres,
hyperbolic spaces, and euclidean spaces of various dimensions.  Under
Ricci flow, the product structure is preserved, and the pieces evolve
in characteristic ways: the spheres each shrink to a point singularity
in finite time (type 1 singularity); the hyperbolic spaces expand for
all time, with no singularity developing; and the euclidean spaces are
flat and static.

Let $g_{S^{n}}$ be the metric on $n$-dimensional sphere $S^n$ with
sectional curvature $1$ and let $g_{H^{n}}$ be the metric on
hyperbolic space $H^n$ with sectional curvature $-1$.  In this section
we again use the notations stated at the beginning of section 2.

\subsection*{B1. $H^{3} \times \mathbb{R}$}

In this case, any initial metric can be written as
\[
g_0=R^2 g_{H^{3}} + du^2 
\] 
for some $R>0$. The Ricci flow solution $g$ is given by
\[
g(t)=(R^{2} + 4 t ) g_{H^{3}} +  du^{2}
\qquad  - \frac{R^{2}}{4} < t < +\infty.
\]

\subsection*{B2. $S^{2} \times \mathbb{R}^{2}$}

In this case,  any initial metric can be written as
\[
g_0=R^2 g_{S^{2}} +  du_1^2 +  du_2^2
\] 
for some $R>0$. The Ricci flow solution $g$ is given by
\[
g(t)=(R^{2} - 2 t ) g_{S^{2}} + du_1^{2} +du_2^2
\qquad  -\infty < t < \frac{R^{2}}{2}.
\]

\subsection*{B3. $H^{2}\times \mathbb{R}^{2}$}

In this case,  any initial metric can be written as
\[
g_0=R^2 g_{H^{2}} +  du_1^2 +  du_2^2
\] 
for some $R>0$. The Ricci flow solution $g$ is given by
\[
g(t)=(R^{2} + 2 t ) g_{H^{2}} +  du_1^{2} +du_2^2
\qquad  - \frac{R^{2}}{2}< t < +\infty .
\]

\subsection*{B4. $S^{2} \times S^{2}$}

In this case,  any initial metric can be written as
\[
g_0= R_1^2 g_{S^{2}}(x) + R_2^2 g_{S^{2}}(y) 
\] 
for some $R_1>0$ and $R_2>0$. The Ricci flow solution
$g$ is given by
\[
g(t)=(R_1^{2} -2t) g_{S^{2}}(x) +  (R_2^{2} - 2 t ) g_{S^{2}} (y)
\qquad -\infty < t < \min \{ \frac{R_1^{2}}{2}, 
\frac{R_2^{2}}{2} \} .
\]

\subsection*{B5. $S^{2} \times H^{2}$}

In this case,  any initial metric can be written as
\[
g_0= R_1^2 g_{S^{2}} + R_2^2 g_{H^{2}} 
\] 
for some $R_1>0$ and $R_2>0$. The Ricci flow solution $g$ is given by
\[
g(t)=(R_1^{2} -2t) g_{S^{2}} +  (R_2^{2} + 2 t ) g_{H^{2}} 
\qquad -\frac{R_2^{2}}{2} < t <  \frac{R_1^{2}}{2}.
\]

\subsection*{B6. $H^{2} \times H^{2}$}

In this case,  any initial metric can be written as
\[
g_0= R_1^2 g_{H^{2}}(x) + R_2^2 g_{H^{2}}(y)
\] 
for some $R_1>0$ and $R_2>0$. The Ricci flow solution
$g$ is given by
\[
g(t)=(R_1^{2} +2t) g_{H^{2}}(x) +  (R_2^{2} + 2 t ) g_{H^{2}}(y)
\qquad \max \{- \frac{R_1^{2}}{2}, 
-\frac{R_2^{2}}{2} \} < t < +\infty .
\]

\subsection*{B7. $\mathbb{C}P^{2}$}

Let $g_{FS}$ be the Fubini-Study metric on
$\mathbb{C}P^{2}$
with constant holomorphic bisectional curvature $+1$. 
Then the Ricci curvature 
$R_{i\bar{j}}(g_{FS})=3(g_{FS})_{i\bar{j}}$.
In this case, any initial metric can be written as (see \cite{KN},
p.277)
\[
g_0= R^2 g_{FS}
\] 
for some $R>0$. The Ricci flow solution $g$ (not the volume-normalized
K\"{a}hler Ricci flow) is given by
\[
g(t)=(R^{2} -6t) g_{FS}
\qquad  -\infty < t < \frac{R^{2}}{6}.
\]

Note that K\"{a}hler Ricci flow with positive holomorphic bisectional
curvature on $\mathbb{C}P^{2}$ has been studied by Chen and Tian
(\cite{CT}); they prove that the (volume-normalized) K\"{a}hler Ricci
flow converges exponentially fast to a K\"{a}hler metric of constant
holomorphic bisectional curvature.
\subsection*{B8. $\mathbb{C}H^{2}$}

Let $g_{\mathbb{C}H^{2}}$ be the K\"{a}hler metric on
$\mathbb{C}H^{2}$
with constant holomorphic bisectional curvature $-1$. 
Then the Ricci curvature 
$R_{i\bar{j}}(g_{\mathbb{C}H^{2}})=
-3(g_{\mathbb{C}H^{2}})_{i\bar{j}}$.
In this case, any initial metric can be written as (see \cite{KN},
p.277)
\[
g_0= R^2 g_{\mathbb{C}H^{2}}
\] 
for some $R>0$. The Ricci flow solution $g$ (not the volume normalized
K\"{a}hler Ricci flow) is given by
\[
g(t)=(R^{2} +6t) g_{\mathbb{C}H^{2}}
\qquad  -\frac{R^{2}}{6} < t < +\infty .
\]

\subsection*{B9. $S^{4}$}

In this case, any initial metric can be written as
\[
g_0= R^2 g_{S^{4}} 
\] 
for some $R>0$. The Ricci flow solution $g$ is given by
\[
g(t)=(R^{2} -6t) g_{S^{4}}  \qquad
 -\infty < t < \frac{R^{2}}{6}   .
\]

\subsection*{B10. $H^{4}$}

In this case, any initial metric can be written as
\[
g_0= R^2 g_{H^{4}} 
\] 
for some $R>0$. The Ricci flow solution $g$ is given by
\[
g(t)=(R^{2} +6t) g_{H^{4}} \qquad
 - \frac{R^{2}}{6} < t < + \infty.
\]


\section*{4. Conclusion}

We have analyzed the Ricci flow for compact four dimensional
homogeneous geometries for which an initial diagonal metric remains
diagonal under the flow.  We obtain explicit solutions in most cases.
We find that if the solution has long-time existence, then it is a
Type III singularity solution.  For volume-normalized flow, there are
examples of collapse to dimensions 1, 2, and 3.

For the nondiagonal cases, the relevant ordinary differential equation
systems are of a similar nature but considerably more complicated.
Numerical techniques should be useful for verifying if the behavior is
similar to that of the diagonal cases.

\bigskip 
\noindent 
\textbf{Acknowledgements}\\
This work was supported in part by National Science Foundation grants
PHY--0354659 and DMS--0405255 at the University of Oregon.  P.L.
thanks McKenzie Wang for some helpful discussions.



\begin{thebibliography}{99999}

\bibitem[B]{B}A. Besse, \emph{Einstein manifolds},
Springer-Verlag, Berlin, 1987.

\bibitem[CT]{CT}X.X. Chen and G. Tian, \emph{
Ricci flow on K\"{a}hler-Einstein surfaces}, Invent. Math. 
\textbf{147} (2002), 487--544.

\bibitem[Ha86]{Ha86}R. Hamilton, \emph{Four-manifolds with positive
curvature operator}, J. Diff. Geom. \textbf{24} (1986), 153--179.

\bibitem[Ha95]{Ha95}R. Hamilton, \emph{The formation of singularities
in the Ricci flow}. Surveys in differential geometry, Vol. II
(Cambridge,
MA, 1993), 7-136, Internat. Press, Cambridge, MA, 1995.

\bibitem[Ha97]{Ha97}R. Hamilton, \emph{Four-manifolds with positive
isotropic curvature}, Comm. Anal. Geom. \textbf{5} (1997), 1-92.

\bibitem[Hi]{Hi}J.A. Hillman, \emph{Four-manifolds, geometries and
knots},
Geometry \& Topology Publications, Coventry, 2002.

\bibitem [Hu]{Hu2}G. Huisken, \emph{Ricci deformation of the metric
on a
Riemannian manifold}, J. Diff. Geom. \textbf{21} (1985), 47-62.

\bibitem[IJ]{IJ}J. Isenberg and M. Jackson, \emph{The Ricci flow of
locally
homogeneous geometries of closed manifolds}, J. Diff. Geom.
\textbf{35}
(1992), 723-741.

\bibitem[I]{I}S. Ishihara, \emph{Homogeneous spaces of four
dimensions}, J. 
Math. Soc. Japan \textbf{7} (1955), 151-168.

\bibitem[K]{K}D. Knopf, \emph{Quasi-convergence of the Ricci flow},
Comm.
Anal. Geom. \textbf{2} (2000), 375-391.

\bibitem[KM]{KM}D. Knopf and K. McLeod, \emph{Quasi-convergence of
model
geometries under the Ricci flow}, Comm. Anal. Geom. \textbf{9}
(2001), 879-919.

\bibitem[KN]{KN}S. Kobayashi and K. Nomizu, \emph{Foundations of
differential geometry}, Vol. II. John Wiley \& Sons, Inc., New York,
1969.

\bibitem[M]{M}M. MacCallum, \emph{On the classification of the real
four-dimensional Lie algebra}, in On Einstein's path, ed. by A.
Harvey,
Springer (1992), 299-317.

\bibitem[Mi]{Mi}J. Milnor, \emph{Curvatures of left invariant metrics
on Lie
groups}, Adv. Math. \textbf{21} (1976), 293-329.

\bibitem[P1]{P1}G. Perelman, \emph{The entropy formula for the Ricci
flow and
its geometric applications}, arXiv:math/0211159.

\bibitem[P2]{P2}G. Perelman, \emph{Ricci flow with surgery on
three-manifolds}
, arXiv:math/0303109.

\bibitem[S]{S}P. Scott, \emph{The geometries of 3-manifolds}, Bull.
London
Math. Soc. \textbf{15} (1983), 401-487.

\bibitem[W]{W}C.T.C. Wall, \emph{Geometries and geometrc structures
in real dimension
 4 and complex dimension 2}, Lecture Notes in Math. \textbf{1167},
401-487, Springer,
 Berlin, 1985.
\end{thebibliography}
\end{document}